\documentclass[12pt, reqno]{amsart}
\usepackage[margin=0.78in]{geometry}
\usepackage[pagewise]{lineno}

\usepackage[english]{babel}
\usepackage[alphabetic, initials]{amsrefs}
\usepackage{amsmath,amssymb,amsfonts,amsthm,enumerate}
\usepackage{url}
\usepackage{graphicx,epstopdf,color}

\numberwithin{equation}{section}

\newtheorem{theorem}{Theorem}[section]

\newcommand{\R}{\mathbb{R}}

\newcommand{\HH}{\mathcal{H}^{n-1}}

\renewcommand{\epsilon}{\varepsilon}

\newcommand{\e}{\varepsilon}

\renewcommand{\leq}{\leqslant}
\renewcommand{\le}{\leqslant}

\renewcommand{\ge}{\geqslant}

\makeatletter
\@namedef{subjclassname@2020}{%
  \textup{2020} Mathematics Subject Classification}
\makeatother

\title[Some nonlocal formulas inspired by an identity of James Simons]{Some nonlocal formulas \\ inspired by an identity of James Simons}
\thanks{DS is supported by the
Australian Future Fellowship
FT230100333.
JT is supported by an Australian Government Research Training Program Scholarship. EV is supported by the Australian Laureate Fellowship FL190100081 {\em Minimal surfaces, free
boundaries, and partial differential equations}. It is a pleasure to thank Joaquim Serra
for many inspiring discussions
and Adam Thompson,
{\em magister geometriae},
for his advice on matters related to differential geometry.}

\author{Serena Dipierro}
\address[Serena Dipierro] {Department of Mathematics and Statistics, The University of Western Australia, 35 Stirling Highway, Crawley, Perth, WA 6009, Australia}
\email{serena.dipierro@uwa.edu.au}

\author{Jack Thompson}
\address[Jack Thompson] {Department of Mathematics and Statistics, The University of Western Australia, 35 Stirling Highway, Crawley, Perth, WA 6009, Australia}
\email{jack.thompson@research.uwa.edu.au}

\author{Enrico Valdinoci}
\address[Enrico Valdinoci]{Department of Mathematics and Statistics, The University of Western Australia, 35 Stirling Highway, Crawley, Perth, WA 6009, Australia}
\email{enrico.valdinoci@uwa.edu.au}

\begin{document}

\begin{abstract}
Inspired by a classical identity proved by James Simons, we establish a new geometric formula
in a nonlocal, possibly fractional, setting.

Our formula also recovers the classical case in the limit, thus providing an approach to Simons' work
that does not heavily rely on differential geometry.
\end{abstract}

\subjclass[2020]{53A10, 49Q05, 35R11.}
\keywords{Geometric identities, minimal surfaces, nonlocal minimal surfaces.}

\maketitle

\tableofcontents

\section{Introduction}

\subsection{Taking inspiration from Simons' work}
A classical identity proved by James Simons in~\cite{MR0233295}
states that at every point of a smooth hypersurface in Euclidean space there holds \begin{align}
    \Delta h &= \operatorname{Hess} H+ H h^2-c^2h .\label{ORIGINAL SIM}
\end{align} Here \(h\) and \(H\) are the second fundamental form and mean curvature with respect to a local choice of unit normal, and, in coordinates, \(h^2\) is the symmetric 2-tensor given by \((h^2)_{ij}=\sum_{k,\ell=1}^{n-1}g^{k\ell}h_{ik}h_{\ell j}\) and \(c^2=\sum_{k,\ell=1}^{n-1} g^{k\ell }(h^2)_{k\ell} \) is the norm squared of the second fundamental form. In particular, it follows from~\eqref{ORIGINAL SIM}, that for a smooth hypersurface with vanishing mean curvature, we have that
\begin{equation*}
\Delta c^2+2c^4=2\sum_{i,j,k=1}^{n-1} |\delta_{k}h_{ij}|^2.
\end{equation*}
Here above, $\delta_{k}$ denotes the tangential derivative
in the $k$-th coordinate direction,
$h_{ij}$ the entries of the second
fundamental form, $c$ the norm of the second fundamental form, and~$\Delta $
the Laplace-Beltrami operator on the hypersurface (see e.g.
formula~(2.16) in~\cite{MR2780140}, or, equivalently, the seventh formula in display on page~123
of~\cite{MR775682}, for further details on this classical formula).

Simons's Identity is pivotal, since it provides the essential ingredient to establish the regularity of stable minimal surfaces up to dimension~$7$.

In this note we speculate about possible generalizations of Simons' Identity to nonlocal settings.
In particular, we will consider the case of boundary of sets
and of level sets of functions. These cases are motivated, respectively, by the study of nonlocal minimal surfaces and nonlocal phase transition equations. The prototypical case of these problems comes from fractional minimal surfaces, as introduced in~\cite{MR2675483}, and we recall that the full regularity theory of the minimizers of the fractional perimeter is one of the main open problems in the field of nonlocal equations:
up to now, this regularity is only known when the ambient space has dimension~$2$, see~\cite{MR3090533}, or up to dimension~$7$
provided that the fractional parameter is sufficiently close to integer, see~\cite{MR3107529}, or when the surface possesses a graphical structure, see~\cite{MR3934589} (see also~\cite{MR3981295, MR4116635, HAR}
for the case of stable nonlocal minimal surfaces, i.e. for surfaces of vanishing nonlocal mean curvature with
nonnegative definite second variation of the corresponding energy functional).

The problem of nonlocal minimal surfaces can also be considered for more general kernels than the one of purely fractional type, see~\cite{MR3930619, MR3981295}, and it can be recovered as the limit in the $\Gamma$-convergence sense of long-range phase coexistence problems, see~\cite{MR2948285}.
In this regard, the regularity properties of nonlocal minimal surfaces are intimately related to the flatness of nonlocal phase transitions, which is also a problem of utmost importance in the contemporary research:
up to now, these flatness properties have been established in dimension up to~$3$, or up to~$8$ for mildly nonlocal operators
under an additional limit assumption, or in dimension~$4$ for the square root of the Laplace operator, see~\cite{MR2177165, MR2498561, MR2644786, MR3148114, MR3280032, MR3740395, MR3812860, MR3939768, MR4050103, MR4124116}, 
the other cases being widely open.\medskip

In this paper, we will not specifically address these regularity and rigidity problems, but rather focus
on a geometric formula which is closely related to Simons' Identity in the nonlocal scenarios.
The application of this formula for the regularity theory appears to be highly nontrivial, since
careful estimates for the reminder terms are needed (in dimension~$3$, a reminder estimate has been recently put forth in~\cite{HAR}).

An interesting by-product of the formula that we present here is that it recovers the classical Simons' Identity
as a limit case. 
Therefore, our nonlocal formula also provides a new approach towards
the original
Simons' Identity, with a new proof which makes only very limited use
of Riemannian geometry and relies instead on some clever use of the integration by parts.\medskip

We stress that,
from our standpoint, having a different proof of the identity by James Simons based on integral, rather than differential, geometry is,  conceptually important. Namely, for us, the matter is not much to find a ``shorter'' proof, in fact some classical proofs of the identity may be shorter than ours, yet they have to introduce a ``suitable notion of connection'' and leverage the full apparatus of differential geometry (in a sense, they literally stand on the shoulders of giants), while our methods are fully elementary and do not use anything more than very basic calculus.

Moreover, being able to recover important differential geometry properties out of integral calculations is, in our opinion, a culturally relevant fact in itself.

Additionally, in the classical case the identity of Simons can be seen as a special situation of Bochner's formula; however, in the nonlocal scenario the Bochner's formula becomes more involuted and involves objects, such as
a norm of the nonlocal Hessian and a nonlocal Ricci curvature, which do not have a clear geometric interpretation.

We also recall that the classical identity of Simons is actually valid for all minimal submanifolds of arbitrary codimension, but a full understanding of nonlocal minimal surfaces is currently lacking
in the general lower dimensional setting (see however~\cite{MR3733825, MR4661775} for suitable notions of nonlocal curvatures).

On top of that, it is highly desirable to find new methods to attack classical problems in differential geometry using novel ideas and unconventional perspectives: in this sense, nonlocal equations seem to provide an excellent training camp for truly new approaches and they often reveal special features distinguishing the classical from the nonlocal worlds, see in particular~\cite{HAR, JSEMA}
(however, to attack the nonlocal case, one has to specifically address situations lacking a ``suitable notion of connection''
and a ``user-friendly Bochner's formula'').
\medskip

Let us now dive into the technical details of our results.

\subsection{The geometric case}\label{PRE1}

Let~$K$ be a radially symmetric kernel satisfying
\begin{equation}\label{ipotesi}
\begin{split}
& K\in C^1(\R^n\setminus\{0\}),\\
& |K(x)|\le \frac{C}{|x|^{n+s}}\\
{\mbox{and }}\qquad &{ |\nabla K(x)|\le\frac{C}{|x|^{n+s+1}}},
\end{split}
\end{equation}
for some~$C>0$ and~$s\in(0,1)$.

Given a set~$E$ with smooth boundary, we consider the $K$-mean curvature of~$E$
at~$x\in\partial E$ given by
\begin{equation}\label{defHs}
H_{K,E}(x):=\frac12\,
\int_{\R^n} \big(\chi_{\R^n\setminus E}(y)-\chi_E(y)\big)\,K(x-y)\,dy.\end{equation}
Notice that the above integral is taken in the principal value sense.

The classical mean curvature of~$E$
will be denoted by~$H_E$.
We define
\begin{equation}\label{defc}
{c_{K,E}}(x):=\sqrt{ \frac12\,\int_{\partial E} \big( \nu_{E}(x)-
\nu_{E}(y)\big)^2\,K(x-y)\,d\HH_y },\end{equation}
being~$\nu_E=(\nu_{E,1},\dots,\nu_{E,n})$ the exterior unit normal of~$E$. The quantity~${c_{K,E}}$
plays in our setting the role played by the norm of the second
fundamental form in the classical case, and we can consider it
the $K$-total curvature of~$E$.

We also define the (minus) $K$-Laplace-Beltrami operator along~$\partial E$
of a function~$f$ by
\begin{equation}\label{defop}
{L_{K,E}} f(x):=\int_{\partial E} \big(f(x)-f(y)\big)\,K(x-y)\,d\HH_y.\end{equation}

As customary, we consider the tangential derivative
\begin{equation}\label{DEL} 
\delta_{E,i} f(x):= \partial_i f(x)-\nu_{E,i}(x)\,\nabla f(x)\cdot \nu_{E}(x)\end{equation}
and we recall that
\begin{equation}\label{SYMMECOR}\delta_{E,i}\nu_{E,j}=\delta_{E,j}\nu_{E,i},\end{equation}
see e.g.
formula~(10.11) in~\cite{MR775682}.

In this setting, our nonlocal formula inspired by Simons' Identity goes as follows:

\begin{theorem}\label{simons}
Let~$K$ be as in~\eqref{ipotesi}.
Let~$E\subset\R^n$ be a set
with smooth boundary and~$x\in\partial E$ with~$\nu_{E}(x)=(0,\dots,0,1)$.

Assume that there exist~$R_0>0$ and~$\beta\in[0, n+s)$ such that
for all~$R\ge R_0$ it holds that
\begin{equation}\label{DEA:PA:0}
\int_{\partial E\cap B_R(x)} \big(|H_E(y)|+1\big)\,d{\mathcal{H}}^{n-1}_y\le CR^\beta,
\end{equation}
for some~$C>0$.

Then, for any~$i$, $j\in\{1,\dots,n-1\}$
it holds that
\begin{equation}\label{SIM:FORNEW}\begin{split}
\delta_{E,i} \delta_{E,j} H_{K,E}(x) \,=\,& -{L_{K,E}} \delta_{E,j} \nu_{E,i} (x) +
c^2_{K,E}(x)\, \delta_{E,j} \nu_{E,i}(x) \\&
-\int_{\partial E} \Big( H_E(y)K(x-y)  -\nu_{E}(y) \cdot \nabla K(x-y) \Big)
\nu_{E,i}(y)\, \nu_{E,j}(y)\,d\HH_y .\end{split}
\end{equation}
\end{theorem}

The proof of Theorem~\ref{simons} will be given in detail in Section~\ref{PROOFSI}.

It is interesting to remark that the result
of Theorem~\ref{simons} ``passes to the limit efficiently
and localizes'': for instance, if one takes~$\rho\in C^\infty_0([-1,1])$,
$\e>0$ and a kernel of the form~$K_\e(x):=\e^{-n-2}\rho(|x|/\e)$,
then, using Theorem~\ref{simons} and sending~$\e\searrow0$,
one recovers the classical Simons' Identity
in~\cite{MR0233295} (such passage to the limit can be performed
e.g. with the analysis in~\cite{MR3230079} and Appendix~C in~\cite{LAWSON}).

The details\footnote{We also remark that condition~\eqref{DEA:PA:0} is obviously
satisfied with~$\beta: =n-1$ when the set~$E$ is smooth and bounded.
For minimizers, and, more generally, stable critical points, of the nonlocal perimeter functional, one still has perimeter estimates (see formula~(1.16)
in Corollary~1.8 of~\cite{MR3981295}): however, in this general case,
estimating the mean curvature, or, in greater generality, the ``second derivatives'' of the set, may be a demanding task, see~\cite{HAR} for some results in this direction.} on how to reconstruct the classical Simons' Identity in the appropriate
limit are given in Section~\ref{APPENDICE A}.

\subsection{Back to the original Simons' Identity}\label{APPENDICE A}

As mentioned above, our nonlocal
formula~\eqref{SIM:FORNEW}
in Theorem~\ref{simons} recovers, in the limit, the original
Simons' Identity proved in~\cite{MR0233295}. The precise result goes as follows:

\begin{theorem}\label{ILTKC}
Let~$E\subset\R^n$ have a smooth boundary.

Let~$x\in\partial E$ and assume that there exist~$R_0>0$ and~$\beta\in[0, n+1)$ such that
for all~$R\ge R_0$ it holds that
\begin{equation}\label{DEA:PA} {\mathcal{H}}^{n-1}\big(
\partial E\cap B_R(x) \big)\le CR^\beta,
\end{equation}
for some~$C>0$. Then, the identity in~\eqref{ORIGINAL SIM}
holds true as a consequence of formula~\eqref{SIM:FORNEW}.
\end{theorem}

The proof of Theorem~\ref{ILTKC} is contained in Section~\ref{OJSDLN-wf}.

We point out that Theorems~\ref{simons} and~\ref{ILTKC} also provide
a new proof of the original Simons' Identity.
Remarkably, our proof relies less on the differential geometry structure
of the hypersurface and it is, in a sense, ``more extrinsic'':
these facts allow us to exploit similar methods also
for the case of 
integrodifferential equations, as will be done in the forthcoming
Section~\ref{IDE2}.

\subsection{The case of integrodifferential equations}\label{IDE2}

The framework that we provide here is
a suitable modification of that given in Section~\ref{PRE1}
for sets. The idea is to ``substitute'' the volume measure~$\chi_E(x)\,dx$
with~$u(x)\,dx$ and the area measure~$\chi_{\partial E}(x)d\HH_x$
with~$|\nabla u(x)|\,dx$. However, one cannot really exploit
the setting of Section~\ref{PRE1} as it is
also for integrodifferential equations, and it is necessary
to ``redo the computation'', so to extrapolate the correct
operators and stability conditions for the solutions.

The technical details go as follows. 
Though more general cases can be considered, for the sake of concreteness, we focus on a kernel~$K$ satisfying
\begin{equation}\label{ipotesi2}
\begin{split}
& K\in C^1(\R^n)\cap L^1(\R^n),\\
&|\nabla K|\in L^1(\R^n),\\
{\mbox{and}}\qquad & K(x)=K(-x).
\end{split}
\end{equation}

Given a function~$u\in W^{1,\infty}(\R^n)$ whose level sets~$\{u=t\}$ are
smooth for a.e.~$t\in\R$,
we define the $K$-mean curvature of~$u$
at~$x\in\R^n$ by
\begin{equation}\label{defHs:2-NN}
H_{K,u}(x):= C_K-
\int_{\R^n} u(y)\,K(x-y)\,dy,\qquad{\mbox{where }}\;
C_K:=\frac12\,\int_{\R^n} K(y)\,dy.
\end{equation}
The setting in~\eqref{defHs:2-NN} has to be compared with~\eqref{defHs}
and especially with the forthcoming formula~\eqref{defHs:2}. The classical mean curvature
of the level sets of~$u$ will be denoted by~$H_u$
(i.e., if~$t_x:=u(x)$,
then~$H_{u}(x)$ is the classical mean curvature of the set~$\{ u>t_x\}$
at~$x$).

We also define the the $K$-total curvature of~$u$ as
\begin{equation}\label{defc-NN}
{c_{K,u}}(x):=\sqrt{ \frac12\,\int_{\R^n} \big( \nu_{u}(x)-
\nu_{u}(y)\big)^2\,K(x-y)\,d\mu_{u,y} },\end{equation}
being~$\nu_{u}(x)$ the exterior unit normal of the level set of~$u$
passing through~$x$ (i.e., if~$t_x:=u(x)$,
then~$\nu_{u}(x)$ is the exterior normal of the set~$\{ u>t_x\}$
at~$x$). In~\eqref{defc-NN}, we also used the notation
\begin{equation}\label{MU}
d\mu_{u,y} :=|\nabla u(y)|\,dy.\end{equation}
Of course, the definition in~\eqref{defc-NN} has to be compared with
that in~\eqref{defc}.
Moreover, by construction we have that
\begin{equation} \label{NUFUP}\nu_{u}(x)=-\frac{\nabla u(x)}{|\nabla u(x)|},\end{equation}
the minus sign coming from the fact that the external derivative
of~$\{u>t_x\}$ points towards points with ``decreasing values'' of~$u$.

We also define the $K$-Laplace-Beltrami operator induced by~$u$
acting on a function~$f$ by
\begin{equation}\label{defop-NN}
{L_{K,u}} f(x):=\int_{\R^n} \big(f(x)-f(y)\big)\,K(x-y)\,d\mu_{u,y}.\end{equation} 
Once again, one can compare~\eqref{defop}
and~\eqref{defop-NN}. Also, we denote by~$\delta_{u,i}$ the tangential
derivatives along the level sets of~$u$ (recall~\eqref{DEL}).
This setting turns out to be
the appropriate one to translate Theorem~\ref{simons}
into a result for solutions of integrodifferential
equations, as will be presented in the forthcoming result:

\begin{theorem}\label{simons-NN}
Let~$K$ be as in~\eqref{ipotesi2}.
Let~$u\in W^{1,\infty}(\R^n)$ and assume
that~$\{u=t\}$ is a
smooth hypersurface with bounded mean curvature for a.e.~$t\in\R$.
For any~$x\in\R^n$ with~$\nu_{u}(x)=(0,\dots,0,1)$
and any~$i$, $j\in\{1,\dots,n-1\}$,
it holds that
\begin{equation}\label{duji483765bv9875bv76c98uoixro}\begin{split}
\delta_{u,i} \delta_{u,j} H_{K,u}(x) \,=\,& -{L_{K,u}} \delta_{u,j} \nu_{u,i} (x) +
c_{K,u}^2(x)\, \delta_{u,j} \nu_{u,i}(x) \\&
-\int_{\R^n} \Big( H_u(y)K(x-y)  -\nu_{u}(y) \cdot \nabla K(x-y) \Big)
\nu_{u,j}(y)\, \nu_{u,i}(y)\,d\mu_{u,y} .\end{split}\end{equation}
\end{theorem}

The proof of Theorem~\ref{simons-NN} is a careful variation
of that of Theorem~\ref{simons}, but, for the sake
of clarity, we provide full details in Section~\ref{ALSMsd}.
We also observe that the choice~$u:=\chi_E$ would formally allow one to recover
Theorem~\ref{simons} from Theorem~\ref{simons-NN}.

\subsection{Stable sets}

In the study of variational problems, a special role is played by the ``stable'' critical points, i.e.
those critical points at which the second derivative of the energy functional is nonnegative definite, see e.g.~\cite{MR3838575}.

In this spirit, in the study of nonlocal minimal surfaces we say that~$E$ is a stable set in~$\Omega$ if~$H_{K,E}(x)=0$ for any~$x\in\Omega\cap\partial E$
and
\begin{equation}\label{ojwfe034}
\frac12\int_{\partial E}\int_{\partial E}(f(x)-f(y))^2\, K(x-y)\,d{\mathcal{H}}^{n-1}_x\,d{\mathcal{H}}^{n-1}_y-
\int_{\partial E} c_{K,E}^2(x)\,f^2(x)\,d{\mathcal{H}}^{n-1}_x\ge0\end{equation}
for any~$f\in C^\infty_0(\Omega)$.

In connection with this, we set
$$ B_{K,E}(u,v ; x) := \frac1 2 \int_{\partial E}(u(x)-u(y))(v(x)-v(y))\,K(x-y)\,d \mathcal H^{n-1}_y,$$
where the integral is taken in the principal value sense, and
$$ B_{K,E} (u,v):= \int_{\partial E} B_{K,E}(u,v ; x) \,d \mathcal H^{n-1}_x . $$
In this notation, the first term in~\eqref{ojwfe034} takes the form~$B_K(f,f)$.

Also, we consider the integrodifferential operator \(L_{K,E}\) previously introduced in~\eqref{defop}. When~$K(x)=\frac1{|x|^{n+s}}$, this operator reduces to the fractional Laplacian, up to normalizing constants.

With this notation, we have:
\begin{theorem}\label{Jfqydwvfbe923ejfn}
Let~$E\subset\R^n$ be a set
with smooth boundary with~$H_{K,E}(y)=0$ for each~$y\in\partial E$. 
Assume that~$E$ is stable in~$\R^n$.

Then, for all \(\eta \in C^\infty_0(\partial E)\), \begin{align*}&
-\int_{\partial E} \bigg \{  \frac12 L_{K,E} c_{K,E}^2(x) +B_{K,E}(c_{K,E},c_{K,E};x) - c_{K,E}^4(x) \bigg \} \eta^2(x) \,d \mathcal H^{n-1}_x \\&\qquad\qquad\qquad\leqslant \int_{\partial E} c_{K,E}^2(x) B_{K,E}(\eta,\eta ; x ) \,d \mathcal H^{n-1}_x .
\end{align*}
\end{theorem}

For the classical counterpart of the above in equality, see e.g.~\cite[equation~(19)]{MR3838575}.
\medskip

The rest of this paper contains the proofs of the results stated above.
Before undertaking the details of the proofs, we mention that the idea of recovering classical results in geometry as a limit of fractional ones, thus providing a unified approach between different disciplines, can offer interesting perspectives 
(see also~\cite{s1KTRA} for limit formulas related to trace problems
and~\cite{KAG} for a recovery technique of the Divergence Theorem coming from a nonlocal perspective).

\section{Proof of Theorem~\ref{simons}}\label{PROOFSI}

Up to a translation, we can suppose that~$0\in\partial E$ and prove
Theorem~\ref{simons} at the origin, hence
we can choose coordinates such that
\begin{equation}\label{ZEROZ}
\nu_{E}(0)=(0,\dots,0,1). 
\end{equation}

We point out that assumption~\eqref{DEA:PA:0} guarantees that all the terms in~\eqref{SIM:FORNEW} are finite, see e.g. the forthcoming technical calculation in~\eqref{d9o3285v456746578fdghdshgxcvdheritv}.

Moreover, we take~$K$ to be smooth, compactly supported and nonsingular, so to be able to
take derivatives inside the integral (the general case then
follows\footnote{The gist of this approximation is,
when $K$ is a singular kernel, to pick a small $\delta \in \left(0, \frac{1}{2}\right)$ and an even mollifier $\eta_\delta \in C^\infty(\R, [0,1])$ such that~$\eta_\delta = 1$ on $[0, \delta] \cup \left[\frac{1}{\delta}, +\infty\right)$,
$\eta_\delta = 0$ on $\left[2\delta, \frac{1}{2\delta}\right)$,
and~$\|\eta_\delta'\|_{L^\infty(\mathbb{R})} \le \frac{1}{\delta}$,
with $\eta_\delta$ pointwise decreasing to $0$ as $\delta \searrow 0$.

Then, one considers a ``regularized kernel'' defined by~$
K_\delta(x) := \left(1 - \eta_\delta(|x|)\right) K(x)$.
The technical advantage is that $K_\delta$ removes the singularity at the origin and is compactly supported, making calculations easier and allowing the general results to be obtained in the limit as $\delta \searrow 0$.} by approximation, see e.g.~\cite[Proposition~6.3]{MR3322379}).
In this way, we rewrite~\eqref{defHs} as
\begin{equation}\label{defHs:2}
H_{K,E}(x)= C_K-
\int_{E} K(x-y)\,dy,\qquad{\mbox{where }}\;
C_K:=\frac12\,\int_{\R^n} K(y)\,dy.
\end{equation}
Also, this is a good definition for all~$x\in\R^n$ (and not only for~$x\in\partial E$),
so we can consider the full gradient of such an expression.
Moreover, for a fixed~$x\in\R^n$, we use the notation
\begin{equation}\label{NOTA1}
\phi(y):=K(x-y).
\end{equation}
In this way, we have that, for any~$\ell\in\{1,\dots,n\}$,
\begin{equation}\label{NOTA2} \partial_\ell K(x-y)=-\partial_\ell\phi(y).\end{equation}
Exploiting this,~\eqref{defHs:2} and the Gauss-Green
Theorem, we see that, for any~$\ell\in\{1,\dots,n\}$,
\begin{equation*}
\begin{split}
& \partial_\ell H_{K,E}(x)=-\int_{E} \partial_\ell K(x-y)\,dy
=\int_{E} \partial_\ell \phi(y)\,dy
=\int_{E} {\rm div}\big(\phi(y) e_\ell\big)\,dy\\
&\qquad
=\int_{\partial E} \nu_{E}(y)\cdot\big(\phi(y) e_\ell\big)\,d\HH_y
=\int_{\partial E} \nu_{E,\ell}(y)\,K(x-y)\,d\HH_y.
\end{split}\end{equation*}
This gives that, for any~$x\in\partial E$,
\begin{equation}\label{UNO}
\nabla H_{K,E}(x)= \int_{\partial E} \nu_{E}(y)\,K(x-y)\,d\HH_y.
\end{equation}
In addition, from~\eqref{defc},
\begin{equation}\label{UNOprimo}
\begin{split}
c^2_{K,E}(x)\,&=\frac12\,\int_{\partial E} \big( \nu_{E}(x)-\nu_{E}(y)\big)^2\,K(x-y)\,d\HH_y\\&
= \int_{\partial E} K(x-y)\,d\HH_y-
\nu_{E}(x)\cdot\int_{\partial E}\nu_{E}(y)\,K(x-y)\,d\HH_y.\end{split}
\end{equation}

Now, we fix the indices~$i$, $j\in\{1,\dots,n-1\}$
and we make use of~\eqref{DEL} and~\eqref{UNO} to find that
\begin{equation}\label{DUE}
\begin{split}
\delta_{E,i} H_{K,E}(x)&=
\partial_i H_{K,E}(x)-\nu_{E,i}(x)\,\nabla H_{K,E}(x)\cdot \nu_{E}(x)\\
&=\int_{\partial E} \nu_{E,i}(y)\,K(x-y)\,d\HH_y
-\nu_{E,i}(x)\,\nu_{E}(x)\cdot\int_{\partial E} \nu_{E}(y)\,K(x-y)\,d\HH_y.
\end{split}\end{equation}
We take another tangential derivative of~\eqref{DUE}
and evaluate it at the origin, recalling~\eqref{ZEROZ}
(which, in particular, gives that~$\nu_{E,i}(0)=0=\nu_{E,j}(0)$
for any~$i$, $j\in\{ 1,\dots,n-1\}$).
In this way, recalling~\eqref{DEL}, we obtain that
\begin{equation}\label{TRE}
\begin{split}&
\delta_{E,j}\delta_{E,i} H_{K,E}(0)\\=\;&
\partial_j\delta_{E,i} H_{K,E}(0)
\\ =\;& \partial_j \left[
\int_{\partial E} \nu_{E,i}(y)\,K(x-y)\,d\HH_y
-\nu_{E,i}(x)\,\nu_{E}(x)\cdot\int_{\partial E} \nu_{E}(y)\,K(x-y)\,d\HH_y\right]_{x=0}
\\ =\;&
\int_{\partial E} \nu_{E,i}(y)\,\partial_j K(-y)\,d\HH_y
-\partial_j\nu_{E,i}(0)\,\nu_{E}(0)\cdot\int_{\partial E} \nu_{E}(y)\,K(-y)\,d\HH_y.
\end{split}\end{equation}
Also, using the notation in~\eqref{NOTA1} and~\eqref{NOTA2}
with~$x:=0$ and~\eqref{DEL}, we see that
\begin{equation}\label{XTRE}
\begin{split} &
\int_{\partial E} \nu_{E,i}(y)\,\partial_j K(-y)\,d\HH_y\\ 
=\;&-\int_{\partial E} \nu_{E,i}(y)\,\partial_j \phi(y)\,d\HH_y\\
=\;&-\int_{\partial E} \nu_{E,i}(y)\,\delta_{E,j} \phi(y)\,d\HH_y-
\int_{\partial E} \nu_{E,i}(y)\,\nu_{E,j}(y)\,\nabla \phi(y)\cdot \nu_{E}(y)\,d\HH_y.
\end{split}\end{equation}

Now we recall an integration by parts formula for tangential derivatives
(see e.g. the first formula\footnote{We stress that the normal on page~122 of~\cite{MR775682}
is internal, according to the distance setting on page~120 therein. This causes
in our notation a sign change with respect to the setting in~\cite{MR775682}.
Also, in the statement of Lemma~10.8 on page~121 in~\cite{MR775682} there is a typo
(missing a mean curvature inside an integral). We also observe that formula~\eqref{DIV:TH}
can also be seen as a version of the Tangential Divergence Theorem, see e.g.
Appendix~A in~\cite{MR2024995}.}
in display on page~122 of~\cite{MR775682}),
namely
\begin{equation}\label{DIV:TH}
\int_{\partial E} \delta_{E,j} f(y)\,d\HH_y=
\int_{\partial E} H_E(y)\,\nu_{E,j}(y)\,f(y)\,d\HH_y,
\end{equation}
being~$H_E$ the classical mean curvature of~$\partial E$.
Applying this formula to the product of two functions, we find that
\begin{equation}\label{EQ:6}
\begin{split}
\int_{\partial E} \delta_{E,j} f(y)\,g(y)\,d\HH_y+
\int_{\partial E} f(y)\,\delta_{E,j} g(y)\,d\HH_y\,&
=
\int_{\partial E} \delta_{E,j} (fg)(y)\,d\HH_y
\\ &=
\int_{\partial E} H_E(y)\,\nu_{E,j}(y)\,f(y)\,g(y)\,d\HH_y.
\end{split}\end{equation}
Using this and~\eqref{NOTA1} (with~$x:=0$ here), we see that
\begin{eqnarray*}&&
-\int_{\partial E} \nu_{E,i}(y)\,\delta_{E,j} \phi(y)\,d\HH_y\\&=&
\int_{\partial E} \delta_{E,j}\nu_{E,i}(y)\,\phi(y)\,d\HH_y
-\int_{\partial E} H_E(y)\,\nu_{E,i}(y)\,\nu_{E,j}(y)\,\phi(y)\,d\HH_y\\
&=&
\int_{\partial E} \delta_{E,j}\nu_{E,i}(y)\,K(-y)\,d\HH_y
-\int_{\partial E} H_E(y)\,\nu_{E,i}(y)\,\nu_{E,j}(y)\,K(-y)\,d\HH_y.\end{eqnarray*}
So, we insert this information into~\eqref{XTRE} and we conclude that
\begin{equation*}
\begin{split}
& \int_{\partial E} \nu_{E,i}(y)\,\partial_j K(-y)\,d\HH_y=
\int_{\partial E} \delta_{E,j}\nu_{E,i}(y)\,K(-y)\,d\HH_y\\&\qquad
-\int_{\partial E} H_E(y)\,\nu_{E,i}(y)\,\nu_{E,j}(y)\,K(-y)\,d\HH_y+
\int_{\partial E} \nu_{E,i}(y)\,\nu_{E,j}(y)\,\nabla K(-y)\cdot \nu_{E}(y)\,d\HH_y.
\end{split}\end{equation*}
Plugging this into~\eqref{TRE}, we get that
\begin{equation}\label{X4}\begin{split}
\delta_{E,j}\delta_{E,i} H_{K,E}(0)=&
\int_{\partial E} \delta_{E,j}\nu_{E,i}(y)\,K(-y)\,d\HH_y
-\int_{\partial E} H_E(y)\,\nu_{E,i}(y)\,\nu_{E,j}(y)\,K(-y)\,d\HH_y\\ &\qquad\quad+
\int_{\partial E} \nu_{E,i}(y)\,\nu_{E,j}(y)\,\nabla K(-y)\cdot \nu_{E}(y)\,d\HH_y\\ &\qquad\quad
-\partial_j\nu_{E,i}(0)\,\nu_{E}(0)\cdot\int_{\partial E} \nu_{E}(y)\,K(-y)\,d\HH_y.
\end{split}\end{equation}

In addition, from~\eqref{UNOprimo},
$$ \partial_j\nu_{E,i}(0)\, c^2_{K,E}(0)
=\int_{\partial E} \partial_j\nu_{E,i}(0)\,K(-y)\,d\HH_y-
\partial_j\nu_{E,i}(0)\,\nu_{E}(0)\cdot\int_{\partial E} \nu_{E}(y)\,K(-y)\,d\HH_y.
$$
Comparing with~\eqref{X4}, we conclude that
\begin{equation*}\begin{split}
\delta_{E,j}\delta_{E,i} H_{K,E}(0)=&
\int_{\partial E} \Big(\delta_{E,j}\nu_{E,i}(y)-\delta_{E,j}\nu_{E,i}(0)\Big)\,K(-y)\,d\HH_y\\ &\qquad
-\int_{\partial E} H_E(y)\,\nu_{E,i}(y)\,\nu_{E,j}(y)\,K(-y)\,d\HH_y\\ &\qquad+
\int_{\partial E} \nu_{E,i}(y)\,\nu_{E,j}(y)\,\nabla K(-y)\cdot \nu_{E}(y)\,d\HH_y
+\partial_j\nu_{E,i}(0)\, c^2_{K,E}(0).
\end{split}\end{equation*}
{F}rom this identity and the definition in~\eqref{defop}, the desired
result plainly follows.\hfill$\Box$

\section{Proof of Theorem~\ref{simons-NN}}\label{ALSMsd}

The proof is similar to that of Theorem~\ref{simons}. Full details are provided for the reader's facility.
Up to a translation, we can prove
Theorem~\ref{simons-NN} at the origin and suppose that
\begin{equation}\label{ZEROZ-NN}
\nu_{u}(0)=(0,\dots,0,1). 
\end{equation}
We observe that our assumptions on the kernel in~\eqref{ipotesi2} yield that all
the terms\footnote{In particular, notice that, by~\eqref{ipotesi2}
and the assumptions
that~$u$ is Lipschitz continuous and the mean curvature of its level sets are bounded,
$$\left|\int_{\R^n} H_u(y)K(x-y)\nu_{u,j}(y)\, \nu_{u,i}(y)\,d\mu_{u,y}\right|\le
\int_{\R^n} \big| H_u(y)K(x-y) \big|\,|\nabla u(y)|\,dy<+\infty.$$} in~\eqref{duji483765bv9875bv76c98uoixro} are finite.

Using~\eqref{defHs:2-NN},
\eqref{MU} and~\eqref{NUFUP}, we see that, for any~$x\in\R^n$,
\begin{equation}\label{UNO-NN}
\begin{split}
&\nabla H_{K,u}(x)= \nabla\left(C_K-
\int_{\R^n} u(x-y)\,K(y)\,dy\right)
=-\int_{\R^n} \nabla u(x-y)\,K(y)\,dy
\\ &\qquad=-\int_{\R^n} \nabla u(y)\,K(x-y)\,dy
=\int_{\R^n} \nu_{u}(y)\,K(x-y)\,d\mu_{u,y}.
\end{split}\end{equation}
In addition, from~\eqref{defc-NN},
\begin{equation}\label{UNOprimo-NN}
\begin{split}
c^2_{K,u}(x)\,&=
\frac12\,\int_{\R^n} \big( \nu_{u}(x)-\nu_{u}(y)\big)^2\,K(x-y)\,d\mu_{u,y}
\\&
= 
\int_{\R^n} K(x-y)\,d\mu_{u,y}
-\nu_{u}(x)\cdot\int_{\R^n} \nu_{u}(y)\,K(x-y)\,d\mu_{u,y}.\end{split}
\end{equation}
Also, in view of~\eqref{DEL} and~\eqref{UNO-NN},
\begin{equation}\label{DUE-NN}
\begin{split}
\delta_{u,i} H_{K,u}(x)&=
\partial_i H_{K,u}(x)-\nu_{u,i}(x)\,\nabla H_{K,u}(x)\cdot \nu_{u}(x)\\
&=\int_{\R^n} \nu_{u,i}(y)\,K(x-y)\,d\mu_{u,y}
-\nu_{u,i}(x)\,\nu_{u}(x)\cdot\int_{\R^n} \nu_{u}(y)\,K(x-y)\,d\mu_{u,y}.
\end{split}\end{equation}
Consequently, using~\eqref{ZEROZ-NN} and~\eqref{DUE-NN}, for all~$i$, $j\in\{1,\dots,n-1\}$,
\begin{equation}\label{TRE-NN}
\begin{split}&
\delta_{u,j}\delta_{u,i} H_{K,u}(0)\\
=\;&
\partial_j\delta_{u,i} H_{K,u}(0)
\\ =\;& \partial_j \left[
\int_{\R^n} \nu_{u,i}(y)\,K(x-y)\,d\mu_{u,y}
-\nu_{u,i}(x)\,\nu_{u}(x)\cdot\int_{\R^n} \nu_{u}(y)\,K(x-y)\,d\mu_{u,y}
\right]_{x=0}
\\ =\;& 
\int_{\R^n} \nu_{u,i}(y)\,\partial_j K(-y)\,d\mu_{u,y}
-\partial_j\nu_{u,i}(0)\,\nu_{u}(0)\cdot\int_{\R^n} \nu_{u}(y)\,K(-y)\,d\mu_{u,y}
.\end{split}\end{equation}

Now, recalling
the notation in~\eqref{NOTA1} and~\eqref{NOTA2}
with~$x:=0$ and~\eqref{DEL}, we obtain that
\begin{equation}\label{XTRE-NN}
\begin{split}&
\int_{\R^n} \nu_{u,i}(y)\,\partial_j K(-y)\,d\mu_{u,y}\\
=\;&-\int_{\R^n} \nu_{u,i}(y)\,\partial_j \phi(y)\,d\mu_{u,y}
\\
=\;&-\int_{\R^n} \nu_{u,i}(y)\,\delta_{u,j} \phi(y)\,d\mu_{u,y}-
\int_{\R^n} \nu_{u,i}(y)\,\nu_{u,j}(y)\,\nabla \phi(y)\cdot \nu_{u}(y)\,d\mu_{u,y}.
\end{split}\end{equation}
Furthermore,
exploiting the Coarea Formula twice
and the tangential integration
by parts identity in~\eqref{EQ:6}, we obtain that
\begin{eqnarray*}
&&
-\int_{\R^n} \nu_{u,i}(y)\,\delta_{u,j} \phi(y)\,d\mu_{u,y}
= -\int_{\R^n} |\nabla u(y)|\,\nu_{u,i}(y)\,\delta_{u,j} \phi(y)\,dy\\
&&\qquad=
-\int_{\R} \int_{\{u(y)=t\}} 
\nu_{u,i}(y)\,\delta_{u,j} \phi(y)\,d\HH_y\,dt
\\ &&\qquad=
\int_{\R} \int_{\{u(y)=t\}} 
\delta_{u,j}\nu_{u,i}(y)\, \phi(y)\,d\HH_y \,dt\\&&\qquad\qquad\qquad
-
\int_{\R} \int_{\{u(y)=t\}} 
H_u(y)\,\nu_{u,i}(y)\,\nu_{u,j}(y)\,\phi(y)\,d\HH_y \,dt
\\ &&\qquad=
\int_{\R^n}|\nabla u(y)|\,
\delta_{u,j}\nu_{u,i}(y)\, \phi(y)\,dy
-\int_{\R^n}|\nabla u(y)|\,
H_u(y)\,\nu_{u,i}(y)\,\nu_{u,j}(y)\,\phi(y)\,dy
\\ &&\qquad=
\int_{\R^n}
\delta_{u,j}\nu_{u,i}(y)\, \phi(y)\,d\mu_{u,y}
-\int_{\R^n}H_u(y)\,\nu_{u,i}(y)\,\nu_{u,j}(y)\,\phi(y)\,d\mu_{u,y}
\\ &&\qquad=
\int_{\R^n}
\delta_{u,j}\nu_{u,i}(y)\, K(-y)\,d\mu_{u,y}
-\int_{\R^n}H_u(y)\,\nu_{u,i}(y)\,\nu_{u,j}(y)\,K(-y)\,d\mu_{u,y}.\end{eqnarray*}
We can now insert this identity into~\eqref{XTRE-NN} and we get that
\begin{equation*}
\begin{split}
& \int_{\R^n} \nu_{u,i}(y)\,\partial_j K(-y)\,d\mu_{u,y}
=
\int_{\R^n} \delta_{u,j}\nu_{u,i}(y)\,K(-y)\,d\mu_{u,y}\\&\qquad
-\int_{\R^n} H_u(y)\,\nu_{u,i}(y)\,\nu_{u,j}(y)\,K(-y)\,d\mu_{u,y}+
\int_{\R^n} \nu_{u,i}(y)\,\nu_{u,j}(y)\,\nabla K(-y)\cdot \nu_{u}(y)\,d\mu_{u,y}.
\end{split}\end{equation*}
Plugging this into~\eqref{TRE-NN} we get that
\begin{equation}\label{X4-NN}\begin{split}&
\delta_{u,j}\delta_{u,i} H_{K,u}(0)=
\int_{\R^n} \delta_{u,j}\nu_{u,i}(y)\,K(-y)\,d\mu_{u,y}
-\int_{\R^n} H_u(y)\,\nu_{u,i}(y)\,\nu_{u,j}(y)\,K(-y)\,d\mu_{u,y}\\ &\qquad\quad+
\int_{\R^n} \nu_{u,i}(y)\,\nu_{u,j}(y)\,\nabla K(-y)\cdot \nu_{u}(y)\,d\mu_{u,y}
-\partial_j\nu_{u,i}(0)\,\nu_{u}(0)\cdot\int_{\R^n} \nu_{u}(y)\,K(-y)\,d\mu_{u,y}.
\end{split}\end{equation}

Also, from~\eqref{UNOprimo-NN}, we have that
$$ \partial_j\nu_{u,i}(0)\, c^2_{K,u}(0)
=\int_{\R^n} \partial_j\nu_{u,i}(0)\,K(-y)\,d\mu_{u,y}-
\partial_j\nu_{u,i}(0)\,\nu_{u}(0)\cdot\int_{\R^n} \nu_{u}(y)\,K(-y)\,d\mu_{u,y}.
$$
Hence, from this and~\eqref{X4-NN}, we conclude that
\begin{equation*}\begin{split}&
\delta_{u,j}\delta_{u,i} H_{K,u}(0)=
\int_{\R^n} \Big(\delta_{u,j}\nu_{u,i}(y)-\delta_{u,j}\nu_{u,i}(0)\Big)\,K(-y)
\,d\mu_{u,y}
-\int_{\R^n} H_u(y)\,\nu_{u,i}(y)\,\nu_{u,j}(y)\,K(-y)\,d\mu_{u,y}
\\ &\qquad\quad+
\int_{\R^n} \nu_{u,i}(y)\,\nu_{u,j}(y)\,\nabla K(-y)\cdot \nu_{u}(y)\,d\mu_{u,y}
+\partial_j\nu_{u,i}(0)\, c^2_{K,u}(0).
\end{split}\end{equation*}
This and~\eqref{defop-NN} give the desired
result.~\hfill$\Box$

\section{Proof of Theorem~\ref{ILTKC}}\label{OJSDLN-wf}
For clarity, we denote by~$\Delta_{\partial E}$
the Laplace-Beltrami operator on the hypersurface~$\partial E$,
by~$\delta_{k,E}$ the tangential derivative
in the $k$th coordinate direction, by~$\nu_E$ the external derivative and by~$c_E$ the norm of the second fundamental form.

To obtain~\eqref{ORIGINAL SIM} as a limit of~\eqref{SIM:FORNEW},
we focus on a special kernel. Namely,
given~$\e>0$, we let
\begin{equation}\label{KERNTOSIM} K_\e(y):=\frac{\e}{|y|^{n+1-\e}}.\end{equation}

We now recall a simple, explicit calculation:
\begin{eqnarray}
&& \label{0191:0}\int_{B_1} x_1^4\,dx=\frac{3\, {\mathcal{H}}^{n-1}(S^{n-1}) }{n(n+2)(n+4) }\\ {\mbox{and }}\,
\label{0191:1} &&\int_{S^{n-1}} \vartheta_1^4 \,d{\mathcal{H}}^{n-1}_\vartheta= 
\frac{3\, {\mathcal{H}}^{n-1}(S^{n-1}) }{n(n+2) } .
\end{eqnarray}
Not to interrupt the flow of the arguments, we postpone 
the proof of formulas~\eqref{0191:0} and~\eqref{0191:1} to Appendix~\ref{jxiewytdenfrhgfndbvhfdnbvfdefgbhwektheru}.

To complete the proof of Theorem~\ref{ILTKC}, without loss of generality, we assume that~$0=x\in\partial E$
and that~$\partial E\cap B_{r_0}$ is the graph of a function~$f:\R^{n-1}\to\R$
with vertical normal, hence~$f(0)=0$ and~$\partial_i f(0)=0$
for all~$i\in\{1,\dots,n-1\}$.
We can also diagonalize the Hessian matrix of~$f$ at~$0$, and obtain that
the mean curvature~$H_E$ at the origin coincides with the trace\footnote{We stress
that we are not dividing the quantity in~\eqref{DEF H COS} by~$n-1$,
to be consistent with the notation in formula~(10.12)
in~\cite{MR775682}.}
of such matrix,
namely
\begin{equation}\label{DEF H COS} H_E(0)= -\big(
\partial^2_1 f(0)+\dots+\partial^2_{n-1} f(0) \big).\end{equation}
The sign convention here is inferred by the assumption that~$E$ is locally the subgraph
of~$f$ and the normal is taken to point outwards.
Consequently, for every~$y=(y',f(y'))\in\partial E\cap B_{r_0}$,
\begin{equation}\label{0-120-A5bis}\begin{split}
&f(y')=\frac12\sum_{i=1}^{n-1} \partial^2_i f(0)\,y_i^2+
O(|y'|^3),\\
&\nabla f(y')=(\partial_1 f(y'),\dots,
\partial_{n-1} f(y'))=\big( \partial^2_1 f(0)\,y_1,\dots,
\partial^2_{n-1} f(0)\,y_{n-1}\big)+
O(|y'|^2)\end{split}\end{equation}
and
\begin{equation}\label{SVIL-NU}
\begin{split}&
\nu_E(y) =\frac{(-\nabla f(y'), 1)}{\sqrt{1+|\nabla f(y')|^2}}=(-\nabla f(y'), 1)+O(|y'|^2)
\\&\qquad=
\big( -\partial^2_1 f(0)\,y_1,\dots,
-\partial^2_{n-1} f(0)\,y_{n-1},1\big)+
O(|y'|^2)\big).\end{split}
\end{equation}
Here, the notation \(g = O(h(\vert y' \vert))\) means that \(\vert g \vert \le 
C \vert h(\vert y' \vert) \vert\) for \(\vert y' \vert \) sufficient close to \(0\) with \(C\) independent of \(\varepsilon\), that is, \(g\) is \emph{uniformly in }\(\varepsilon\) big O of \(h\) as \(\vert y'\vert \to 0\).  As a consequence,
for every~$y=(y',f(y'))\in\partial E\cap B_{r_0}$,
\begin{equation} \label{SVIPRO}|y|^2=|y'|^2+|f(y')|^2=|y'|^2+O(|y'|^4)=
|y'|^2\big(1+O(|y'|^2)\big),\end{equation}
and, for any~$i,j\in\{1,\dots,n-1\}$,
\begin{equation} \label{DOPPIO NU}\nu_{E,j}(y)\nu_{E,i}(y)=
\partial^2_{j} f(0)\partial^2_{i} f(0) y_j y_i
+
O(|y'|^3) .\end{equation}
Thus, using~\eqref{SVIPRO}, we see that, for any fixed~$\alpha\in\R$,
\begin{equation} \label{SVIPRO1}|y|^\alpha=
|y'|^\alpha\big(1+O(|y'|^2)\big)^{\alpha/2}=
|y'|^\alpha\big(1+O(|y'|^2)\big).\end{equation}

Then, from~\eqref{SVIL-NU} and~\eqref{SVIPRO1},
we obtain that, for any~$\ell\in\{1,\dots,n-1\}$ and~$y\in\partial E\cap B_{r_0}$,
\begin{eqnarray*}
\nu_{E,\ell}(y)\partial_\ell K_\e(-y)&=&-
\frac{(n+1-\e)\e \partial^2_\ell f(0)\,y_\ell^2}{|y|^{n+3-\e}}+\e O\left( |y'|^{\e-n}\right)\\&=&-
\frac{(n+1-\e)\e \partial^2_\ell f(0)\,y_\ell^2}{|y'|^{n+3-\e}}+\e O\left( |y'|^{\e-n}\right)
\end{eqnarray*}
and also, recalling \eqref{0-120-A5bis},
\begin{eqnarray*}
\nu_{E,n}(y)\partial_n K_\e(-y)&=&
\frac{(n+1-\e)\e y_n(1+O(|y'|^2))}{|y|^{n+3-\e}}\\&=&\frac12\sum_{\ell=1}^{n-1}
\frac{(n+1-\e)\e \partial^2_\ell f(0)\,y_\ell^2}{|y'|^{n+3-\e}}
+\e O(|y'|^{\e-n}).
\end{eqnarray*}
Accordingly, we have that
$$ \nu_{E}(y) \cdot \nabla K_\e(-y) =-\frac12\sum_{\ell=1}^{n-1}
\frac{(n+1-\e)\e \partial^2_\ell f(0)\,y_\ell^2}{|y'|^{n+3-\e}}
+\e O(|y'|^{\e-n}).$$ 
We thereby deduce from the latter identity and~\eqref{DOPPIO NU}
(and exploiting an odd symmetry argument)
that, for any~$r\in(0,r_0]$,
\begin{equation}\label{LAPRIINT}\begin{split}&
-\frac2{(n+1-\e)\e}
\int_{\partial E\cap B_r} \nu_{E}(y) \cdot \nabla K_\e(-y) 
\nu_{E,i}(y)\, \nu_{E,j}(y)\,d\HH_y\\ =\;&
\sum_{\ell=1}^{n-1} \int_{\partial E\cap B_r}\left(
\frac{ \partial^2_\ell f(0)
\,\partial^2_j f(0)\,\partial^2_i f(0)
\,y_\ell^2\,y_j\,y_i}{|y'|^{n+3-\e}}+O(|y'|^{2+\e-n})\right)
\,d\HH_y\\
=\;&
\sum_{\ell=1}^{n-1} \int_{\{|y'|<r\}} \left(
\frac{ \partial^2_\ell f(0)
\,\partial^2_j f(0)\,\partial^2_i f(0)
\,y_\ell^2\,y_j\,y_i}{|y'|^{n+3-\e}}+O(|y'|^{2+\e-n})\right)
\,\sqrt{1+|\nabla f(y')|^2}\,dy'\\
=\;&
\sum_{\ell=1}^{n-1} \int_{\{|y'|<r\}} \left(
\frac{ \partial^2_\ell f(0)
\,\partial^2_j f(0)\,\partial^2_i f(0)
\,y_\ell^2\,y_j\,y_i}{|y'|^{n+3-\e}}+O(|y'|^{2+\e-n})\right)
\,dy'\\
=\;&
\sum_{\ell=1}^{n-1} \int_{\{|y'|<r\}} \left(
\frac{ \partial^2_\ell f(0)
\,(\partial^2_j f(0))^2
\,y_\ell^2\,y_j^2\,\delta_{ji}}{|y'|^{n+3-\e}}+O(|y'|^{2+\e-n})\right)
\,dy'.
\end{split}\end{equation}

Furthermore, exploiting again~\eqref{DOPPIO NU}
and~\eqref{SVIPRO1}, we see that
\begin{equation}\label{LAPRIINT2}
\begin{split}&
\int_{\partial E\cap B_r} H_E(y)K_\e(-y)\,\nu_{E,i}(y)\, \nu_{E,j}(y)\,d\HH_y\\=\;&
\e\int_{\partial E\cap B_r} \left(\frac{H_E(y)\partial^2_{j} f(0)\partial^2_{i} f(0) y_j y_i}{
|y'|^{n+1-\e}}+O(|y'|^{2+\e-n})\right)\,d\HH_y
\\=\;&
\e\int_{\{|y'|<r\}}\left( \frac{H_E(0)\partial^2_{j} f(0)\partial^2_{i} f(0) y_j y_i}{
|y'|^{n+1-\e}}+O(|y'|^{2+\e-n})\right)\,dy'
\\=\;&
\e\int_{\{|y'|<r\}} \left(\frac{H_E(0)\,(\partial^2_{j} f(0))^2 \,y_j^2\,\delta_{ji}}{
|y'|^{n+1-\e}}+O(|y'|^{2+\e-n})\right)\,dy'.
\end{split}
\end{equation}
Now we use polar coordinates in~$\R^{n-1}$ to observe that
\begin{equation} \label{RESTO}
\int_{\{|y'|<r\}} |y'|^{2+\e-n}\,dy'=
{\mathcal{H}}^{n-2}(S^{n-2})\,\int_0^r \rho^{2+\e-n}\,\rho^{n-2}\,d\rho
=\frac{C\,r^{1+\e}}{1+\e},\end{equation}
for some~$C>0$.

Moreover, for any fixed index~$j\in\{1,\dots,n-1\}$,
\begin{equation}\label{BAL}\begin{split}&
\e\int_{\{|y'|<r\}} \frac{y_j^2}{|y'|^{n+1-\e}}\,dy'
=\frac{\e}{n-1}\sum_{k=1}^{n-1}\int_{\{|y'|<r\}} \frac{y_k^2}{|y'|^{n+1-\e}}\,dy'\\&\qquad
=\frac{\e}{n-1}\int_{\{|y'|<r\}} \frac{dy'}{|y'|^{n-1-\e}}
=\frac{\e\,{\mathcal{H}}^{n-2}(S^{n-2})}{n-1}\,\int_0^r
\rho^{\e-1}\,d\rho
=\varpi\,r^\e,\end{split}
\end{equation}
where
\begin{equation}\label{VARPI}
\varpi:=\frac{{\mathcal{H}}^{n-2}(S^{n-2})}{n-1}.
\end{equation}

Now, we compute the term~$
\e\int_{\{|y'|<r\}} \frac{y_\ell^2\,y_j^2}{|y'|^{n+3-\e}}\,dy'$.
For this, first of all we deal with the case~$\ell=j$: in this situation, we have that
\begin{equation}\label{891-304959}
\e\int_{\{|y'|<r\}} \frac{y_j^4}{|y'|^{n+3-\e}}\,dy' =
\e\int_{\{|y'|<r\}} \frac{y_1^4}{|y'|^{n+3-\e}}\,dy'=C_\star\,r^\e,
\end{equation}
where
\begin{equation}\label{0239494i292110202}\begin{split}
&C_\star:= \e\int_{\{|y'|<1\}} \frac{y_1^4}{|y'|^{n+3-\e}}\,dy'
= \e\iint_{(\rho,\vartheta)\in(0,1)\times S^{n-2}}
\rho^{\e-1}\,\vartheta_1^4\,d\rho\,d{\mathcal{H}}^{n-2}_\vartheta\\
&\qquad\qquad= \int_{\vartheta\in S^{n-2}}
\vartheta_1^4 \,d{\mathcal{H}}^{n-2}_\vartheta
=
\frac{3\, {\mathcal{H}}^{n-2}(S^{n-2}) }{(n-1)(n+1) }=
\frac{3\varpi }{n+1}
,\end{split}\end{equation}
thanks to~\eqref{0191:1} (applied here in one dimension less).

Moreover,
the number of different indices~$k$, $m\in\{1,\dots,n-1\}$ is equal to~$(n-1)(n-2)$
and so,
for each~$j\ne\ell\in\{1,\dots,n-1\}$,
\begin{equation*}\begin{split}&
\e\int_{\{|y'|<r\}} \frac{y_\ell^2\,y_j^2}{|y'|^{n+3-\e}}\,dy'
=\e\int_{\{|y'|<r\}} \frac{y_1^2\,y_2^2}{|y'|^{n+3-\e}}\,dy'
=\frac{\e}{(n-1)(n-2)}\sum_{k\neq m=1}^{n-1}\int_{\{|y'|<r\}}
\frac{y_k^2\,y_m^2}{|y'|^{n+3-\e}}\,dy'
\\&\qquad=\frac{\e}{(n-1)(n-2)}\left[ \sum_{k, m=1}^{n-1}\int_{\{|y'|<r\}}
\frac{y_k^2\,y_m^2}{|y'|^{n+3-\e}}\,dy' -
\sum_{k=1}^{n-1}\int_{\{|y'|<r\}}
\frac{y_k^4}{|y'|^{n+3-\e}}\,dy'
\right]
\\&\qquad
=\frac{\e}{(n-1)(n-2)}\int_{\{|y'|<r\}} \frac{dy'}{|y'|^{n-1-\e}}
-\frac{C_\star\,r^\e}{n-2}=\frac{\varpi\,r^\e}{n-2}-\frac{C_\star\,r^\e}{n-2}
=\frac{\varpi\,r^\e}{n+1} .\end{split}
\end{equation*}
{F}rom this and~\eqref{891-304959},
we obtain that
\begin{equation*} \e\int_{\{|y'|<r\}} \frac{y_\ell^2\,y_j^2}{|y'|^{n+3-\e}}\,dy'=
\frac{(1+2\delta_{\ell j})\,\varpi\,r^\e}{n+1}.\end{equation*} 
Substituting this identity and~\eqref{RESTO} into~\eqref{LAPRIINT},
and recalling also~\eqref{DEF H COS},
we conclude that
\begin{equation}\label{LAPRIINT:BIS}\begin{split}&-\lim_{\e\searrow0}
\frac2{n+1-\e}
\int_{\partial E\cap B_r} \nu_{E}(y) \cdot \nabla K_\e(-y) 
\nu_{E,i}(y)\, \nu_{E,j}(y)\,d\HH_y\\ =\;&\lim_{\e\searrow0}
\e\sum_{\ell=1}^{n-1} \int_{\{|y'|<r\}} 
\frac{ \partial^2_\ell f(0)
\,(\partial^2_j f(0))^2
\,y_\ell^2\,y_j^2\,\delta_{ji}}{|y'|^{n+3-\e}}
\,dy'
\\ =\;&\lim_{\e\searrow0} \frac{\varpi\,r^\e}{n+1}
\sum_{\ell=1}^{n-1} \partial^2_\ell f(0)
\,(\partial^2_j f(0))^2\,(1+2\delta_{\ell j})\,\delta_{ji}
\\ =\;& -\lim_{\e\searrow0} \frac{\varpi\,r^\e\,H_E(0)}{n+1}
\,(\partial^2_j f(0))^2\,\delta_{ji}
+ \frac{2\varpi\,r^\e}{n+1}\,(\partial^2_j f(0))^3\,\delta_{ji}
\\ =\;& -\frac{\varpi\,H_E(0)}{n+1}
\,(\partial^2_j f(0))^2\,\delta_{ji}
+ \frac{2\varpi}{n+1}\,(\partial^2_j f(0))^3\,\delta_{ji}.
\end{split}\end{equation}
Similarly,
substituting~\eqref{RESTO} and~\eqref{BAL}
into~\eqref{LAPRIINT2}, we obtain that, as~$\e\searrow0$,
\begin{equation*}
\begin{split}&\lim_{\e\searrow0} 
\int_{\partial E\cap B_r} H_E(y)K_\e(-y)\,\nu_{E,i}(y)\, \nu_{E,j}(y)\,d\HH_y\\=\;&\lim_{\e\searrow0} 
\e\int_{\{|y'|<r\}} \frac{H_E(0)\,(\partial^2_{j} f(0))^2 \,y_j^2\,\delta_{ji}}{
|y'|^{n+1-\e}}\,dy'\\=\;&\lim_{\e\searrow0} 
\varpi\,r^\e\,H_E(0)\,(\partial^2_{j} f(0))^2 \,\delta_{ji} \\=\;&
\varpi\,H_E(0)\,(\partial^2_{j} f(0))^2 \,\delta_{ji}.
\end{split}
\end{equation*}
{F}rom this and~\eqref{LAPRIINT:BIS} it follows that
\begin{equation}\label{FI:002:PRE}
\begin{split}&
\lim_{\e\searrow0}\int_{\partial E\cap B_r} \Big( H_E(y)K_\e(-y) 
-\nu_{E}(y) \cdot \nabla K_\e(-y) \Big)
\nu_{E,i}(y)\, \nu_{E,j}(y)\,d\HH_y\\
=\;&\frac{\varpi}{2}\,H_E(0)\,(\partial^2_{j} f(0))^2 \,\delta_{ji}+
\varpi\,(\partial^2_{j} f(0))^3 \,\delta_{ji}.
\end{split}\end{equation}

Now we exploit~\eqref{DEA:PA:0} and we see that
\begin{equation}\begin{split}\label{d9o3285v456746578fdghdshgxcvdheritv}
&\left|
\int_{\partial E\setminus B_r} \Big( H_E(y)K_\e(-y) 
-\nu_{E}(y) \cdot \nabla K_\e(-y) \Big)
\nu_{E,i}(y)\, \nu_{E,j}(y)\,d\HH_y\right|\\
\le\;& C\e\int_{\partial E\setminus B_r}\left( \frac{|H_E(y)|}{|y|^{n+1-\e}}+
\frac1{|y|^{n+2-\e}}\right)\,d\HH_y\\
\le\;& C\e\int_{\partial E\setminus B_r}\frac{1}{|y|^{n+1-\e}}
\left(  |H_E(y)|+{r}^{-1}\right)\,d\HH_y\\
\le\;& C\,(1+r^{-1})\,\e\int_{\partial E\setminus B_r} \frac{|H_E(y)|+1}{|y|^{n+1-\e}}\,d\HH_y\\
=\;& C\,(1+r^{-1})\,\e\,\sum_{k=0}^{+\infty}
\int_{\partial E\cap(B_{2^{k+1}r}\setminus B_{2^kr})} \frac{|H_E(y)|+1}{|y|^{n+1-\e}}\,d\HH_y\\
\le\;& \frac{C\,(1+r^{-1})\,\e}{r^{n+1-\e}}\,\sum_{k=0}^{+\infty}
\frac1{2^{k(n+1-\e)}}
\int_{\partial E\cap(B_{2^{k+1}r}\setminus B_{2^kr})} \big(|H_E(y)|+1\big)\,d\HH_y
\\ \le\;& \frac{C\,(1+r^{-1})\,\e}{r^{n+1-\e}}\,\sum_{k=0}^{+\infty}
\frac{(2^{k+1}r)^\beta}{2^{k(n+1-\e)}}
\\ \le\;& \frac{2^\beta\,C\,(1+r^{-1})\,\e}{r^{n+1-\e-\beta}}\,\sum_{k=0}^{+\infty}
\frac{1}{2^{k(n+1-\e-\beta)}}
\\ \le\;& \frac{2^\beta\,C\,(1+r^{-1})\,\e}{r^{n+1-\e-\beta}}\,\sum_{k=0}^{+\infty}
\frac{1}{2^{\frac{k(n+1-\beta)}2}},
\end{split}\end{equation}
up to renaming~$C$ line after line, 
and the last term in~\eqref{d9o3285v456746578fdghdshgxcvdheritv}
is infinitesimal as~$\e\searrow0$.

Consequently,
\[\lim_{\e\searrow0}\int_{\partial E\setminus B_r} \Big( H_E(y)K_\e(-y) 
-\nu_{E}(y) \cdot \nabla K_\e(-y) \Big)
\nu_{E,i}(y)\, \nu_{E,j}(y)\,d\HH_y=0.\]
This and~\eqref{FI:002:PRE} give that
\begin{equation}\label{FI:002}
\begin{split}&
\lim_{\e\searrow0}\int_{\partial E} \Big( H_E(y)K_\e(-y) 
-\nu_{E}(y) \cdot \nabla K_\e(-y) \Big)
\nu_{E,i}(y)\, \nu_{E,j}(y)\,d\HH_y\\
=\;&\frac{\varpi}{2}\,H_E(0)\,(\partial^2_{j} f(0))^2 \,\delta_{ji}+
\varpi\,(\partial^2_{j} f(0))^3 \,\delta_{ji}.
\end{split}\end{equation}

In addition, from Lemma A.2 of \cite{LAWSON},
we have that
\begin{equation}\label{FI:003}
\lim_{\e\searrow0} {L_{K_\e,E}}=-\frac{\varpi}{2}\,\Delta_{\partial E},
\end{equation}
where the notation in~\eqref{VARPI} has been used.
Similarly, from Lemma A.4 of \cite{LAWSON},
\begin{equation}\label{FI:004}
\lim_{\e\searrow0} c_{K_\e,E}^2=\frac\varpi2\,
c_E^2,\end{equation}
being~$c_E$ the norm of the second fundamental form of~$\partial E$.

Therefore, using~\eqref{FI:002}, \eqref{FI:003}
and~\eqref{FI:004},
we obtain that
\begin{equation}\label{PRE:SI}\begin{split}
&
\lim_{\e\searrow0}\Big[- {L_{K_\e,E}} \delta_{E,j} \nu_{E,i} (0)+
c_{K_\e,E}^2(0)\, \delta_{E,j} \nu_{E,i}(0) \\&\qquad-
\int_{\R^n} \Big( H_E(y)K_\e(-y)  -\nu_{E}(y) \cdot \nabla K_\e(-y) \Big)
\nu_{E,i}(y)\, \nu_{E,j}(y)\,d\HH_y \Big]\\
=\;& \frac{\varpi}{2}\,\Delta_{\partial E}\delta_{E,j} \nu_{E,i} (0)
+\frac{\varpi}2\,c_E^2(0)\, \delta_{E,j} \nu_{E,i}(0) -
\frac{\varpi}{2}\,H_E(0)\,(\partial^2_{j} f(0))^2 \,\delta_{ji}-
\varpi\,(\partial^2_{j} f(0))^3 \,\delta_{ji}.
\end{split}\end{equation}

Now, given two functions~$\psi$, $\phi$, we exploit~\eqref{EQ:6}
twice to obtain that
\begin{equation}\label{017:PRE}
\begin{split}&
\int_{\partial E} \delta_{E,i} \delta_{E,j}\psi(x)\,\phi(x)\,d\HH_x\\
=&-\int_{\partial E} \delta_{E,j}\psi(x)\,\delta_{E,i} \phi(x)\,d\HH_x+
\int_{\partial E} H_E(x)\,\nu_{E,i}(x)\,\delta_{E,j}\psi(x)\,\phi(x)\,d\HH_x\\
=&
\int_{\partial E} \psi(x)\,\delta_{E,j} \delta_{E,i} \phi(x)\,d\HH_x-
\int_{\partial E} H_E(x)\,\nu_{E,j}(x)\,\psi(x)\,\delta_{E,i} \phi(x)\,d\HH_x
\\&\qquad+
\int_{\partial E} H_E(x)\,\nu_{E,i}(x)\,\delta_{E,j}\psi(x)\,\phi(x)\,d\HH_x.
\end{split}\end{equation}
On the other hand, applying~\eqref{EQ:6}
once again, we see that
\begin{eqnarray*}
&& 
\int_{\partial E} H_E(x)\,\nu_{E,i}(x)\,\delta_{E,j}\psi(x)\,\phi(x)\,d\HH_x
\\&=&
-\int_{\partial E} \delta_{E,j}\big(H_E(x)\,\nu_{E,i}(x)\,\phi(x)\big)\,\psi(x)\,d\HH_x
+
\int_{\partial E} H_E^2(x)\,\nu_{E,i}(x)\,\nu_{E,j}(x)\,\psi(x)\,\phi(x)\,d\HH_x.
\end{eqnarray*}
Plugging this information into~\eqref{017:PRE}, we find that
\begin{equation}\label{017}
\begin{split}&
\int_{\partial E} \delta_{E,i} \delta_{E,j}\psi(x)\,\phi(x)\,d\HH_x\\
=&
\int_{\partial E} \psi(x)\,\delta_{E,j} \delta_{E,i} \phi(x)\,d\HH_x-
\int_{\partial E} H_E(x)\,\nu_{E,j}(x)\,\psi(x)\,\delta_{E,i} \phi(x)\,d\HH_x
\\&\qquad-\int_{\partial E} \delta_{E,j}\big(H_E(x)\,\nu_{E,i}(x)\,\phi(x)\big)\,\psi(x)\,d\HH_x\\&\qquad
+
\int_{\partial E} H_E^2(x)\,\nu_{E,i}(x)\,\nu_{E,j}(x)\,\psi(x)\,\phi(x)\,d\HH_x.\end{split}\end{equation}

Applying~\eqref{017} (twice, at the beginning
with~$\psi:=H_{K_\e,E}(x)$ and at the end with~$\psi:=H_{E}(x)$)
and considering~$\phi$ as a test function,
the convergence of~$H_{K_\e,E}$ to~$\frac{\varpi\,H_E}{2}$
(see Theorem~12 in~\cite{MR3230079}) gives that
\begin{eqnarray*}&&
-\lim_{\e\searrow0}
\int_{\partial E} \delta_{E,i} \delta_{E,j}H_{K_\e,E}(x)\,\phi(x)\,d\HH_x
\\&=&\lim_{\e\searrow0}\left[ -
\int_{\partial E} H_{K_\e,E}(x)\,\delta_{E,j} \delta_{E,i} \phi(x)\,d\HH_x+
\int_{\partial E} H_E(x)\,\nu_{E,j}(x)\,H_{K_\e,E}(x)\,\delta_{E,i} \phi(x)\,d\HH_x\right.
\\&&\qquad+ 
\int_{\partial E} \delta_{E,j}\big(H_E(x)\,\nu_{E,i}(x)\,\phi(x)\big)\,H_{K_\e,E}(x)\,d\HH_x
\\&&\left.\qquad-
\int_{\partial E} H_E^2(x)\,\nu_{E,i}(x)\,\nu_{E,j}(x)\,H_{K_\e,E}(x)\,\phi(x)\,d\HH_x
\right]
\\&=& \frac{\varpi}{2}\,\left[-
\int_{\partial E} H_{E}(x)\,\delta_{E,j} \delta_{E,i} \phi(x)\,d\HH_x+
\int_{\partial E} H_E(x)\,\nu_{E,j}(x)\,H_{E}(x)\,\delta_{E,i} \phi(x)\,d\HH_x\right.
\\&&\qquad+
\int_{\partial E} \delta_{E,j}\big(H_E(x)\,\nu_{E,i}(x)\,\phi(x)\big)\,H_{E}(x)\,d\HH_x
\\&&\left.\qquad-
\int_{\partial E} H_E^2(x)\,\nu_{E,i}(x)\,\nu_{E,j}(x)\,H_{E}(x)\,\phi(x)\,d\HH_x
\right]
\\&=& -\frac{\varpi}{2}\,
\int_{\partial E} \delta_{E,i} \delta_{E,j}H_{E}(x)\,\phi(x)\,d\HH_x.
\end{eqnarray*}
This says that~$ \delta_{E,i} \delta_{E,j}H_{K_\e,E} $ converges 
to~$\frac{\varpi}{2}\delta_{E,i} \delta_{E,j}H_{E}$ in the distributional sense
as~$\e\searrow0$. 

We also recall that~$H_{K_\e,E} $ is smooth, thanks
to regularity of~$E$ (see e.g.~\cite[equation~(49)]{MR3331523}), and thus the derivatives
of~$H_{K_\e,E} $ are smooth too. Accordingly,
by the Ascoli-Arzel\`a
Theorem, we have
that~$ \delta_{E,i} \delta_{E,j}H_{K_\e,E} $
converges strongly up to a subsequence, and therefore the uniqueness
of the limit gives that~$ \delta_{E,i} \delta_{E,j}H_{K_\e,E} $ converges also
pointwise
to~$\frac{\varpi}{2}\delta_{E,i} \delta_{E,j}H_{E}$.

Combining this with~\eqref{PRE:SI}, we obtain that
\begin{equation}\label{PRE:PUNTO}
\begin{split}
&
\lim_{\e\searrow0}\Big[ \delta_{E,i} \delta_{E,j}H_{K_\e,E}(0)+
{L_{K_\e,E}} \delta_{E,j} \nu_{E,i} (0)-
c_{K_\e,E}^2(0)\, \delta_{E,j} \nu_{E,i}(0) \\&\qquad+
\int_{\R^n} \Big( H_E(y)K_\e(-y)  -\nu_{E}(y) \cdot \nabla K_\e(-y) \Big)
\nu_{E,i}(y)\, \nu_{E,j}(y)\,d\HH_y \Big]\\
=\;& 
\frac{\varpi}{2}\delta_{E,i} \delta_{E,j}H_{E}(0)
-\frac{\varpi}{2}\,\Delta_{\partial E}\delta_{E,j} \nu_{E,i} (0)
-\frac\varpi2\,c_E^2(0)\, \delta_{E,j} \nu_{E,i}(0) \\&\qquad+
\frac{\varpi}{2}\,H_E(0)\,(\partial^2_{j} f(0))^2 \,\delta_{ji}+
\varpi\,(\partial^2_{j} f(0))^3 \,\delta_{ji}.
\end{split}\end{equation}
By formula~\eqref{SIM:FORNEW}, we know that the left hand side
of~\eqref{PRE:PUNTO} is equal to zero.
Therefore,
\begin{equation}\label{D3:00}\delta_{E,i} \delta_{E,j}H_{E}(0)=
\Delta_{\partial E}\delta_{E,j} \nu_{E,i} (0)
+c_E^2(0)\, \delta_{E,j} \nu_{E,i}(0) -H_E(0)\,(\partial^2_{j} f(0))^2 \,\delta_{ji}-2(\partial^2_{j} f(0))^3 \,\delta_{ji}. \end{equation}
 Since \(H_E\) is a scalar function
defined on the hypersurface~$\partial E$, if  \(\nabla_i\) denotes the \(i\)-th covariant derivative then \begin{align}
    \delta_{E,i} \delta_{E,j}H_{E}(0)&= \nabla_i\nabla_j H_E(0) .\label{D3:01}
 \end{align} 

Next, we claim that, with respect to the coordinates \(y' \mapsto (y',f(y'))\), we have \begin{align}
    \Delta_{\partial E}\delta_{E,j} \nu_{E,i} (0) &= \Delta_{\partial E} h_{ij}(0)+ 2 (\partial^2_j f(0))^3 \delta_{ij}\label{TC61biat} 
\end{align}
where \(h_{ij}\) are the components of the second fundamental form. Once~\eqref{TC61biat} has been established, we obtain
from~\eqref{D3:00} and~\eqref{D3:01} that 
\begin{align*}
    \nabla_i \nabla_jH_{E}(0)&=
\Delta_{\partial E}h_{ij}(0)
+c_E^2(0)\, h_{ij}(0) -H_E(0)\,h_{jj}^2(0) \,\delta_{ij}\\&=
\Delta_{\partial E}h_{ij}(0)
+c_E^2(0)\, h_{ij}(0) -H_E(0)(h^2)_{ij}
\end{align*} 
 where \(h^2\) is the symmetric 2-tensor given by \((h^2)_{ij}=\sum_{k,\ell=1}^{n-1}g^{k\ell}h_{ik}h_{\ell j}\). This is one of the general versions
of the classical identity by Simons in~\cite[Theorem~4.2.1]{MR0233295} and coincides with the
concise formulation given in~\eqref{ORIGINAL SIM}.

Now, to establish~\eqref{TC61biat}, we recall that the metric tensor with respect to the graphical coordinates is \begin{align*}
    g_{ij}(y') = \delta_{ij} + \partial_i f (y')\partial_jf(y'),
\end{align*} so \(g_{ij}(0)=\delta_{ij}\), \(\partial_kg_{ij}(0)=0\), \(\partial_k g^{ij}(0)=0\), and \begin{equation}\begin{split}
        \partial_{k\ell}g_{ij}(0)&= \partial_{ik} f (0)\partial_{j\ell}f(0)+\partial_{i\ell} f (0)\partial_{kj}f(0) \label{QTLrMFcB}\\
        &=(\partial_i^2f(0))(\partial_j^2f(0)) (\delta_{ik}\delta_{j\ell}+\delta_{i\ell}\delta_{kj}).
\end{split}\end{equation} 

Furthermore, if we consider the Christoffel symbols: \begin{align*}
    \Gamma^k_{ij} = \frac 12 \sum_{\alpha=1}^{n-1} g^{\alpha k } \big ( \partial_i g_{\alpha j} + \partial_j g_{i\alpha}-\partial_\alpha g_{ij} \big ) 
\end{align*} then \(\Gamma^k_{ij}(0)=0\) and  \begin{align*}
    \partial_\ell\Gamma^k_{ij}(0)&= \frac 12 \sum_{\alpha=1}^{n-1} g^{\alpha k }(0) \big ( \partial_{i\ell} g_{\alpha j}(0) + \partial_{\ell j} g_{i\alpha}(0)-\partial_{\alpha\ell} g_{ij}(0) \big )\\
    &=\frac 12\big ( \partial_{i\ell} g_{k j}(0) + \partial_{\ell j} g_{ik}(0)-\partial_{k\ell} g_{ij}(0) \big ). 
\end{align*} Hence, it follows from~\eqref{QTLrMFcB} that,
for each~$i$, $j$, $k\in\{1,\dots,n-1\}$,
\begin{align*}
    \partial_k \Gamma^j_{ik}(0)&= \frac 12 \big ( \partial_{k i } g_{jk}(0)+\partial_{kk}g_{ij}(0)-\partial_{kj}g_{ik}(0) \big ) \\
    &= \frac 12 (\partial^2_j f(0))(\partial^2_kf(0)) ( \delta_{ij} \delta_{kk}+\delta_{ik}\delta_{jk})+\frac 12 (\partial^2_i f(0))(\partial^2_jf(0)) ( \delta_{ik} \delta_{jk}+\delta_{ik}\delta_{jk}) \\
    &\qquad - \frac 12  (\partial^2_i f(0))(\partial^2_kf(0)) ( \delta_{ij} \delta_{kk}+\delta_{ik}\delta_{jk}) \\
    &=(\partial^2_j f(0))^2 \delta_{ij}.
\end{align*} Furthermore, the second fundamental form in coordinates is \begin{align*}
    h_{ij}(y')&= - \frac{\partial_{ij}f(y')}{\sqrt{1+\vert \nabla f(y')\vert^2}},
\end{align*} so, in particular, \begin{align*}
    h_{ij}(0)&= - (\partial_j^2f(0)) \delta_{ij}.
\end{align*} Finally, using all of the formulas above and the formula for
covariant derivatives in coordinates (see e.g.\cite[first formula of page~6]{Hebey}), we have at \(0\) that \begin{align}
    \Delta_{\partial E}h_{ij} &= \sum_{\alpha,\beta=1}^{n-1}g^{\alpha \beta } \nabla_\beta (\nabla_\alpha h_{ij}) \nonumber\\
    &= \sum_{\alpha=1}^{n-1}  \partial_\alpha \bigg( \partial_\alpha h_{ij} - \sum_{\sigma=1}^{n-1} \Gamma^\sigma_{i\alpha}h_{\sigma j }-\sum_{\sigma=1}^{n-1}\Gamma^\sigma_{j\alpha}h_{\sigma i } \bigg  )\nonumber \\
    &= \sum_{\alpha=1}^{n-1}\partial_{\alpha \alpha}h_{ij}(0)- \sum_{\alpha,\sigma=1}^{n-1} h_{\sigma j }\partial_\alpha \Gamma^\sigma_{i\alpha}-\sum_{\alpha,\sigma=1}^{n-1}h_{\sigma i } \partial_\alpha\Gamma^\sigma_{j\alpha}\nonumber  \\
    &=  \sum_{\alpha=1}^{n-1}\partial_{\alpha \alpha}h_{ij}(0)+\partial_j^2f(0) \sum_{\alpha=1}^{n-1} \partial_\alpha \Gamma^j_{i\alpha}+\partial_i^2f(0) \sum_{\alpha=1}^{n-1} \partial_\alpha\Gamma^i_{j\alpha}\nonumber \\
    &= \sum_{\alpha=1}^{n-1}\partial_{\alpha \alpha}h_{ij}(0)+2(\partial_i^2f(0))^3 \delta_{ij}. \label{vEgDz27I}
\end{align}

 Next, using~\eqref{SVIL-NU},
 if $1\leq j \leq n-1 $,
  \begin{align*}
    \partial_j \nu_{E,i} (y) &= - \frac{\partial_{ij}f}{\sqrt{1+\vert \nabla f \vert^2}}+ \frac{\partial_if }{(1+\vert \nabla f \vert^2)^{3/2}}\sum_{k=1}^{n-1} \partial_k f \partial_{kj}f, 
\end{align*} and \(\partial_n \nu_{E,i} (y)=0\). 

As a result, for \(1\leq i,j\leq n-1\), \begin{align*}
    \delta_j \nu_{E,i} (y',f(y')) &= - \frac{\partial_{ij}f}{\sqrt{1+\vert \nabla f \vert^2}}+ \frac{\partial_if }{(1+\vert \nabla f \vert^2)^{3/2}}\sum_{k=1}^{n-1} \partial_k f \partial_{kj}f - \sum_{\ell=1}^{n-1}\bigg ( - \frac{\partial_{i\ell}f}{\sqrt{1+\vert \nabla f \vert^2}}\\
    &\qquad + \frac{\partial_if }{(1+\vert \nabla f \vert^2)^{3/2}}\sum_{k=1}^{n-1} \partial_k f \partial_{k\ell}f\bigg ) \frac{\partial_\ell f \partial_jf}{1+\vert \nabla f \vert^2} \\
    &= h_{ij}(y')+E_{ij}(y'),
\end{align*} where \begin{align*}
    E_{ij}(y')&= \frac{1 }{(1+\vert \nabla f \vert^2)^{3/2}}\sum_{k=1}^{n-1} \big ( \partial_if \partial_k f \partial_{kj}f +\partial_kf \partial_{ik}f \partial_kf  \big ) -  \frac{\partial_i f \partial_j f}{(1+\vert \nabla f \vert^2)^{5/2}}\sum_{k,\ell=1}^{n-1}\partial_kf\partial_\ell f \partial_{k\ell }f.
\end{align*} Hence, the metric is orthogonal at \(0\), one has (at \(0\)) \begin{align*}
     \Delta (\delta_j \nu_{E,i}) &= \sum_{\alpha=1}^{n-1 }\partial^2_\alpha (\delta_j \nu_{E,i}(y',f(y')) ) \\
     &= \sum_{\alpha=1}^{n-1 } \big ( \partial^2_\alpha h_{ij}(0)+\partial^2_\alpha E_{ij}(0) \big ) \\
     &= \Delta_{\partial E}h_{ij}(0)-2(\partial_j^2f(0))^3 \delta_{ij}+\sum_{\alpha=1}^{n-1 } \partial^2_\alpha E_{ij}(0) 
\end{align*} and, by a direct calculation \begin{align*}
    \sum_{\alpha=1}^{n-1 } \partial^2_\alpha E_{ij}(0) &= 4(\partial_jf(0))^3 \delta_{ij}   ,
\end{align*} which establishes~\eqref{TC61biat}.~\hfill$\Box$

\section{Proof of Theorem~\ref{Jfqydwvfbe923ejfn}}
Let \(\eta \in C^\infty_0({\partial E})\) be arbitrary. Using \(f: = c_K \eta\) as a test function in~\eqref{ojwfe034}, we have that \begin{eqnarray*}
0&\leqslant& B_K(c_K\eta , c_K \eta ) - \int_{\partial E} c_K^4 \eta^2 \,d \mathcal H^n .
\end{eqnarray*} Next, for all \(x,y\in {\partial E}\), we have that \begin{eqnarray*}
(c_K(x)\eta(x)-c_K(y)\eta(y))^2 &=& \big (c_K(x)(\eta(x)-\eta(y))+ \eta(y) (c_K(x)-c_K(y)) \big )^2\\
&=& c_K^2(x) (\eta(x)-\eta(x))^2 + \eta^2(y)(c_K(x)-c_K(y))^2 \\
&&\qquad + 2c_K(x)\eta(y)(c_K(x)-c_K(y))(\eta(x)-\eta(y)),
\end{eqnarray*}
so it follows that \begin{eqnarray*}
B_K(c_K\eta , c_K \eta ) &=& \int_{\partial E} c_K^2(x) B_K(\eta , \eta ; x) \,d \mathcal H^n_x +  \int_{\partial E} \eta^2(x) B_K(c_K , c_K ; x) \,d \mathcal H^n_x +I  ,
\end{eqnarray*} where \begin{eqnarray*}
I:=  \int_{\partial E} \int_{\partial E} c_K(x)\eta(y)(c_K(x)-c_K(y))(\eta(x)-\eta(y))\,K(x-y)\,d \mathcal H^n_y \,d \mathcal H^n_x.
\end{eqnarray*} Next, by symmetry of \(x\) and \(y\), we have that \begin{eqnarray*}
I= \frac 1 2  \int_{\partial E} \int_{\partial E} (c_K(x)\eta(y)+c_K(y)\eta(x))(c_K(x)-c_K(y))(\eta(x)-\eta(y))\,K(x-y) \,d \mathcal H^n_y \,d \mathcal H^n_x.
\end{eqnarray*}
Moreover, by a simple algebraic manipulation, \begin{eqnarray*}&&
 (c_K(x)\eta(y)+c_K(y)\eta(x))(c_K(x)-c_K(y))(\eta(x)-\eta(y)) \\
&&\qquad= \frac12 (\eta^2(x) -\eta^2(y))(c_K^2(x)-c_K^2(y)) - \frac12 (c_K(x)-c_K(y))^2(\eta(x)-\eta(y))^2 \\
&&\qquad\leqslant \frac12 (\eta^2(x) -\eta^2(y))(c_K^2(x)-c_K^2(y))
\end{eqnarray*}and accordingly \begin{eqnarray*}
I&\leqslant&\frac 1 4  \int_{\partial E} \int_{\partial E} (\eta^2(x) -\eta^2(y))(c_K^2(x)-c_K^2(y))
\,K(x-y)\,d \mathcal H^n_y \,d \mathcal H^n_x \\
&=& \frac 1 2  \int_{\partial E} \int_{\partial E} \eta^2(x)(c_K^2(x)-c_K^2(y))\,K(x-y) \,d \mathcal H^n_y \,d \mathcal H^n_x \\
&=& \frac 12 \int_{\partial E} \eta^2(x) {\mathcal{L}}_K  c_K^2(x) \,d \mathcal H^n_x.
\end{eqnarray*}  Hence, we have that \begin{eqnarray*}&&
B_K(c_K\eta , c_K \eta ) \\&&\qquad\leqslant \int_{\partial E} c_K^2(x) B_K(\eta , \eta ; x) \,d \mathcal H^n_x +  \int_{\partial E} \bigg \{  B_K(c_K , c_K ; x)+\frac 12 \int_{\partial E} \eta^2(x) {\mathcal{L}}_K  c_K^2(x) \bigg \} \eta^2(x)\,d \mathcal H^n_x 
\end{eqnarray*} and the result follows. 
\hfill$\Box$

\appendix

\section{Proof of formulas~\eqref{0191:0} and~\eqref{0191:1}}\label{jxiewytdenfrhgfndbvhfdnbvfdefgbhwektheru}

Let $$
Q:=\int_{B_1} x_1^4\,dx \qquad{\mbox{and}}\qquad D:=\int_{B_1} x_1^2x_2^2\,dx
.$$
We consider the isometry~$x\mapsto X\in \R^n$ given by
$$ X_1:=\frac{x_1-x_2}{\sqrt2},\qquad
X_2:=\frac{x_1+x_2}{\sqrt2},\qquad X_i:=x_i \quad{\mbox{ for all }}\; i\in\{3,\dots,n\}.$$
We notice that
$$ 4X_1^2 X_2^2 = (2X_1X_2)^2=\big( (x_1-x_2)(x_1+x_2)\big)^2=
(x_1^2-x_2^2)^2 = x_1^4+x_2^4-2x_1^2 x_2^2$$
and therefore, by symmetry,
$$ 4D=\int_{B_1} 4X_1^2 X_2^2\,dX=\int_{B_1}
\big( x_1^4+x_2^4-2x_1^2 x_2^2\big)\,dx
=2Q-2D,$$
which gives
\begin{equation}\label{24t-d293}
D=\frac{Q}3.
\end{equation}

On the other hand
$$ |x|^4 = (|x|^2)^2=\left( \sum_{i=1}^n x_i^2\right)^2=
\sum_{i,j=1}^n x_i^2 x_j^2
=\sum_{i=1}^n x_i^4+
\sum_{i\ne j=1}^n x_i^2 x_j^2.$$
Therefore, by polar coordinates and symmetry,
\begin{eqnarray*}
&&\frac{ {\mathcal{H}}^{n-1}(S^{n-1}) }{n+4 }=
{\mathcal{H}}^{n-1}(S^{n-1})
\int_0^1 \rho^{n+3}\,d\rho=
\int_{B_1} |x|^4\,dx\\ &&\qquad=
\sum_{i=1}^n \int_{B_1}x_i^4\,dx+
\sum_{i\ne j=1}^n \int_{B_1} x_i^2 x_j^2\,dx=nQ+n(n-1)D.
\end{eqnarray*}
{F}rom this and~\eqref{24t-d293} we deduce that
$$ \frac{ {\mathcal{H}}^{n-1}(S^{n-1}) }{n+4 } = \frac{n(n+2)\,Q}{3},$$
hence
$$ \frac{3\, {\mathcal{H}}^{n-1}(S^{n-1}) }{n(n+2)(n+4) }=Q
=\int_0^1\int_{S^{n-1}} \rho^{n+3}\vartheta_1^4\,d{\mathcal{H}}^{n-1}_\vartheta\,d\rho=
\frac{1}{n+4}\int_{S^{n-1}} \vartheta_1^4 \,d{\mathcal{H}}^{n-1}_\vartheta,
$$
which gives~\eqref{0191:0} and~\eqref{0191:1}, as desired.

\vfill


\section*{References}
\begin{biblist}


\bib{MR3230079}{article}{
   author={Abatangelo, Nicola},
   author={Valdinoci, Enrico},
   title={A notion of nonlocal curvature},
   journal={Numer. Funct. Anal. Optim.},
   volume={35},
   date={2014},
   number={7-9},
   pages={793--815},
   issn={0163-0563},
   review={\MR{3230079}},
   doi={10.1080/01630563.2014.901837},
}


\bib{MR3331523}{article}{
   author={Barrios, Bego\~{n}a},
   author={Figalli, Alessio},
   author={Valdinoci, Enrico},
   title={Bootstrap regularity for integro-differential operators and its
   application to nonlocal minimal surfaces},
   journal={Ann. Sc. Norm. Super. Pisa Cl. Sci. (5)},
   volume={13},
   date={2014},
   number={3},
   pages={609--639},
   issn={0391-173X},
   review={\MR{3331523}},
}


\bib{MR2644786}{article}{
   author={Cabr\'{e}, Xavier},
   author={Cinti, Eleonora},
   title={Energy estimates and 1-D symmetry for nonlinear equations
   involving the half-Laplacian},
   journal={Discrete Contin. Dyn. Syst.},
   volume={28},
   date={2010},
   number={3},
   pages={1179--1206},
   issn={1078-0947},
   review={\MR{2644786}},
   doi={10.3934/dcds.2010.28.1179},
}


\bib{MR3148114}{article}{
   author={Cabr\'{e}, Xavier},
   author={Cinti, Eleonora},
   title={Sharp energy estimates for nonlinear fractional diffusion
   equations},
   journal={Calc. Var. Partial Differential Equations},
   volume={49},
   date={2014},
   number={1-2},
   pages={233--269},
   issn={0944-2669},
   review={\MR{3148114}},
   doi={10.1007/s00526-012-0580-6},
}


\bib{MR4116635}{article}{
   author={Cabr\'{e}, Xavier},
   author={Cinti, Eleonora},
   author={Serra, Joaquim},
   title={Stable $s$-minimal cones in $\Bbb{R}^3$ are flat for $s\sim 1$},
   journal={J. Reine Angew. Math.},
   volume={764},
   date={2020},
   pages={157--180},
   issn={0075-4102},
   review={\MR{4116635}},
   doi={10.1515/crelle-2019-0005},
}


\bib{MR3934589}{article}{
   author={Cabr\'{e}, Xavier},
   author={Cozzi, Matteo},
   title={A gradient estimate for nonlocal minimal graphs},
   journal={Duke Math. J.},
   volume={168},
   date={2019},
   number={5},
   pages={775--848},
   issn={0012-7094},
   review={\MR{3934589}},
   doi={10.1215/00127094-2018-0052},
}


\bib{MR3838575}{article}{
   author={Cabr\'{e}, Xavier},
   author={Poggesi, Giorgio},
   title={Stable solutions to some elliptic problems: minimal cones, the
   Allen-Cahn equation, and blow-up solutions},
   conference={
      title={Geometry of PDEs and related problems},
   },
   book={
      series={Lecture Notes in Math.},
      volume={2220},
      publisher={Springer, Cham},
   },
   date={2018},
   pages={1--45},
   review={\MR{3838575}},
}


\bib{MR3280032}{article}{
   author={Cabr\'{e}, Xavier},
   author={Sire, Yannick},
   title={Nonlinear equations for fractional Laplacians II: Existence,
   uniqueness, and qualitative properties of solutions},
   journal={Trans. Amer. Math. Soc.},
   volume={367},
   date={2015},
   number={2},
   pages={911--941},
   issn={0002-9947},
   review={\MR{3280032}},
   doi={10.1090/S0002-9947-2014-05906-0},
}


\bib{MR2177165}{article}{
   author={Cabr\'{e}, Xavier},
   author={Sol\`a-Morales, Joan},
   title={Layer solutions in a half-space for boundary reactions},
   journal={Comm. Pure Appl. Math.},
   volume={58},
   date={2005},
   number={12},
   pages={1678--1732},
   issn={0010-3640},
   review={\MR{2177165}},
   doi={10.1002/cpa.20093},
}


\bib{MR2675483}{article}{
   author={Caffarelli, L.},
   author={Roquejoffre, J.-M.},
   author={Savin, O.},
   title={Nonlocal minimal surfaces},
   journal={Comm. Pure Appl. Math.},
   volume={63},
   date={2010},
   number={9},
   pages={1111--1144},
   issn={0010-3640},
   review={\MR{2675483}},
   doi={10.1002/cpa.20331},
}


\bib{MR3107529}{article}{
   author={Caffarelli, Luis},
   author={Valdinoci, Enrico},
   title={Regularity properties of nonlocal minimal surfaces via limiting
   arguments},
   journal={Adv. Math.},
   volume={248},
   date={2013},
   pages={843--871},
   issn={0001-8708},
   review={\MR{3107529}},
   doi={10.1016/j.aim.2013.08.007},
}


\bib{HAR}{article}{
title={Nonlocal approximation of minimal surfaces: optimal estimates from stability},
author={Chan, Hardy},
author={Dipierro, Serena},
author={Serra, Joaquim},
author={Valdinoci, Enrico},
   journal = {arXiv e-prints},
archivePrefix = {arXiv},
       eprint = {2308.06328},
}


\bib{MR3981295}{article}{
   author={Cinti, Eleonora},
   author={Serra, Joaquim},
   author={Valdinoci, Enrico},
   title={Quantitative flatness results and $BV$-estimates for stable
   nonlocal minimal surfaces},
   journal={J. Differential Geom.},
   volume={112},
   date={2019},
   number={3},
   pages={447--504},
   issn={0022-040X},
   review={\MR{3981295}},
   doi={10.4310/jdg/1563242471},
}




\bib{MR2780140}{book}{
   author={Colding, Tobias Holck},
   author={Minicozzi, William P., II},
   title={A course in minimal surfaces},
   series={Graduate Studies in Mathematics},
   volume={121},
   publisher={American Mathematical Society, Providence, RI},
   date={2011},
   pages={xii+313},
   isbn={978-0-8218-5323-8},
   review={\MR{2780140}},
}


\bib{LAWSON}{article}{
   author={D\'{a}vila, Juan},
   author={del Pino, Manuel},
   author={Wei, Juncheng},
   title={Nonlocal $s$-minimal surfaces and Lawson cones},
   journal={J. Differential Geom.},
   volume={109},
   date={2018},
   number={1},
   pages={111--175},
   issn={0022-040X},
   review={\MR{3798717}},
   doi={10.4310/jdg/1525399218},
}


\bib{MR3740395}{article}{
   author={Dipierro, Serena},
   author={Farina, Alberto},
   author={Valdinoci, Enrico},
   title={A three-dimensional symmetry result for a phase transition
   equation in the genuinely nonlocal regime},
   journal={Calc. Var. Partial Differential Equations},
   volume={57},
   date={2018},
   number={1},
   pages={Paper No. 15, 21},
   issn={0944-2669},
   review={\MR{3740395}},
   doi={10.1007/s00526-017-1295-5},
}


	


\bib{MR4124116}{article}{
   author={Dipierro, Serena},
   author={Serra, Joaquim},
   author={Valdinoci, Enrico},
   title={Improvement of flatness for nonlocal phase transitions},
   journal={Amer. J. Math.},
   volume={142},
   date={2020},
   number={4},
   pages={1083--1160},
   issn={0002-9327},
   review={\MR{4124116}},
   doi={10.1353/ajm.2020.0032},
}


\bib{MR2024995}{book}{
   author={Ecker, Klaus},
   title={Regularity theory for mean curvature flow},
   series={Progress in Nonlinear Differential Equations and their
   Applications},
   volume={57},
   publisher={Birkh\"auser Boston, Inc., Boston, MA},
   date={2004},
   pages={xiv+165},
   isbn={0-8176-3243-3},
   review={\MR{2024995}},
}






\bib{MR3322379}{article}{
   author={Figalli, A.},
   author={Fusco, N.},
   author={Maggi, F.},
   author={Millot, V.},
   author={Morini, M.},
   title={Isoperimetry and stability properties of balls with respect to
   nonlocal energies},
   journal={Comm. Math. Phys.},
   volume={336},
   date={2015},
   number={1},
   pages={441--507},
   issn={0010-3616},
   review={\MR{3322379}},
   doi={10.1007/s00220-014-2244-1},
}


\bib{MR4050103}{article}{
   author={Figalli, Alessio},
   author={Serra, Joaquim},
   title={On stable solutions for boundary reactions: a De Giorgi-type
   result in dimension $4+1$},
   journal={Invent. Math.},
   volume={219},
   date={2020},
   number={1},
   pages={153--177},
   issn={0020-9910},
   review={\MR{4050103}},
   doi={10.1007/s00222-019-00904-2},
}
	
\bib{MR775682}{book}{
   author={Giusti, Enrico},
   title={Minimal surfaces and functions of bounded variation},
   series={Monographs in Mathematics},
   volume={80},
   publisher={Birkh\"auser Verlag, Basel},
   date={1984},
   pages={xii+240},
   isbn={0-8176-3153-4},
   review={\MR{775682}},
   doi={10.1007/978-1-4684-9486-0},
}


\bib{s1KTRA}{article}{
 author={Grube, Florian},
   author={Hensiek, Thorben},
   title={Robust nonlocal trace spaces and Neumann problems},
   journal={Nonlinear Anal.},
   volume={241},
   date={2024},
   pages={Paper No. 113481, 35},
   issn={0362-546X},
   review={\MR{4685916}},
   doi={10.1016/j.na.2023.113481},
}



\bib{Hebey}{book}{
title = {Nonlinear analysis on manifolds: Sobolev spaces and inequalities},
	volume = {5},
	isbn = {0-9658703-4-0 0-8218-2700-6},
	series = {Courant Lecture Notes in Mathematics},
	pagetotal = {x+309},
	publisher = {New York University, Courant Institute of Mathematical Sciences, New York; American Mathematical Society, Providence, {RI}},
	author = {Hebey, Emmanuel},
	date = {1999},
	mrnumber = {1688256},
	keywords = {Isoperimetry},
	file = {Hebey - 1999 - Nonlinear analysis on manifolds Sobolev spaces an.pdf:C\:\\Users\\00106759\\Zotero\\storage\\IJXJBJD6\\Hebey - 1999 - Nonlinear analysis on manifolds Sobolev spaces an.pdf:application/pdf},
}

\bib{KAG}{article}{
 author={Hepp, Solveig},
   author={Kassmann, Moritz},
   title={The divergence theorem and nonlocal counterparts},
   journal={Bull. Lond. Math. Soc.},
   volume={56},
   date={2024},
   number={2},
   pages={711--733},
   issn={0024-6093},
   review={\MR{4711580}},
   doi={10.1112/blms.12960},
}


\bib{MR3930619}{book}{
   author={Maz\'{o}n, Jos\'{e} M.},
   author={Rossi, Julio Daniel},
   author={Toledo, J. Juli\'{a}n},
   title={Nonlocal perimeter, curvature and minimal surfaces for measurable
   sets},
   series={Frontiers in Mathematics},
   publisher={Birkh\"{a}user/Springer, Cham},
   date={2019},
   pages={xviii+123},
   isbn={978-3-030-06242-2},
   isbn={978-3-030-06243-9},
   review={\MR{3930619}},
   doi={10.1007/978-3-030-06243-9},
}


\bib{MR3733825}{article}{
   author={Paroni, Roberto},
   author={Podio-Guidugli, Paolo},
   author={Seguin, Brian},
   title={On the nonlocal curvatures of surfaces with or without boundary},
   journal={Commun. Pure Appl. Anal.},
   volume={17},
   date={2018},
   number={2},
   pages={709--727},
   issn={1534-0392},
   review={\MR{3733825}},
   doi={10.3934/cpaa.2018037},
}


\bib{MR4661775}{article}{
   author={Paroni, Roberto},
   author={Podio-Guidugli, Paolo},
   author={Seguin, Brian},
   title={On a notion of nonlocal curvature tensor},
   journal={J. Elasticity},
   volume={154},
   date={2023},
   number={1-4},
   pages={61--79},
   issn={0374-3535},
   review={\MR{4661775}},
   doi={10.1007/s10659-023-09985-w},
}
		
\bib{MR3812860}{article}{
   author={Savin, Ovidiu},
   title={Rigidity of minimizers in nonlocal phase transitions},
   journal={Anal. PDE},
   volume={11},
   date={2018},
   number={8},
   pages={1881--1900},
   issn={2157-5045},
   review={\MR{3812860}},
   doi={10.2140/apde.2018.11.1881},
}


\bib{MR3939768}{article}{
   author={Savin, O.},
   title={Rigidity of minimizers in nonlocal phase transitions II},
   journal={Anal. Theory Appl.},
   volume={35},
   date={2019},
   number={1},
   pages={1--27},
   issn={1672-4070},
   review={\MR{3939768}},
   doi={10.4208/ata.oa-0008},
}


\bib{MR2948285}{article}{
   author={Savin, Ovidiu},
   author={Valdinoci, Enrico},
   title={$\Gamma$-convergence for nonlocal phase transitions},
   journal={Ann. Inst. H. Poincar\'{e} C Anal. Non Lin\'{e}aire},
   volume={29},
   date={2012},
   number={4},
   pages={479--500},
   issn={0294-1449},
   review={\MR{2948285}},
   doi={10.1016/j.anihpc.2012.01.006},
}


\bib{MR3090533}{article}{
   author={Savin, Ovidiu},
   author={Valdinoci, Enrico},
   title={Regularity of nonlocal minimal cones in dimension 2},
   journal={Calc. Var. Partial Differential Equations},
   volume={48},
   date={2013},
   number={1-2},
   pages={33--39},
   issn={0944-2669},
   review={\MR{3090533}},
   doi={10.1007/s00526-012-0539-7},
}


\bib{JSEMA}{article}{
author={Serra, Joaquim},
   title={Nonlocal minimal surfaces: recent developments, applications, and
   future directions},
   journal={SeMA J.},
   volume={81},
   date={2024},
   number={2},
   pages={165--191},
   issn={2254-3902},
   review={\MR{4743530}},
   doi={10.1007/s40324-023-00345-1},
}


\bib{MR0233295}{article}{
   author={Simons, James},
   title={Minimal varieties in riemannian manifolds},
   journal={Ann. of Math. (2)},
   volume={88},
   date={1968},
   pages={62--105},
   issn={0003-486X},
   review={\MR{0233295}},
   doi={10.2307/1970556},
}


\bib{MR2498561}{article}{
   author={Sire, Yannick},
   author={Valdinoci, Enrico},
   title={Fractional Laplacian phase transitions and boundary reactions: a
   geometric inequality and a symmetry result},
   journal={J. Funct. Anal.},
   volume={256},
   date={2009},
   number={6},
   pages={1842--1864},
   issn={0022-1236},
   review={\MR{2498561}},
   doi={10.1016/j.jfa.2009.01.020},
}


\end{biblist}
\end{document}